\documentclass[12pt,reqno]{article}
\usepackage{microtype}
\usepackage{colordvi}
\usepackage{amsmath,amsfonts,amsthm,amssymb,latexsym}
\usepackage{indentfirst}
\allowdisplaybreaks[4]
\usepackage{graphicx}
\usepackage{verbatim}
\usepackage{color}
\usepackage{amscd}
\usepackage{verbatim}
\usepackage{graphicx}

\usepackage{url}

\usepackage{setspace}

\usepackage[left=1.25in, right=1.25in, top=1in, bottom=1in]{geometry}
\usepackage{longtable}

\usepackage{CJKutf8}
\usepackage{appendix}
\usepackage{amsmath,amssymb,amsthm,color,mathrsfs}
\usepackage{enumitem,anysize}%
\usepackage{verbatim}
\usepackage{amssymb}
\usepackage{tikz}
\usepackage{pgfplots}

\marginsize{1in}{1in}{1in}{1in}%

\newtheorem{theorem}{Theorem}[section]
\newtheorem{corollary}[theorem]{Corollary}

\newtheorem{Thm}[theorem]{Theorem}

\newtheorem{lem}[theorem]{Lemma}
\newtheorem{lemma}[theorem]{Lemma}

\newtheorem{example}[theorem]{Example}

\newtheorem{problem}[theorem]{Problem}
\newtheorem{definition}[theorem]{Definition}
\newtheorem{Def}[theorem]{Definition}

\newtheorem{assumption}[theorem]{Assumption}
\newtheorem{assumptions}[theorem]{Assumptions}

\newtheorem{remark}[theorem]{Remark}
\newtheorem{Rem}[theorem]{Remark}
\newtheorem{examples}[theorem]{Examples}
\newtheorem{remarks}[theorem]{Remarks}
\newtheorem{claim*}[theorem]{Claim}
\newtheorem{notation}[theorem]{Notation}

\newcommand{\R}{\mathbb{R}}

\newcommand{\blemma}{\begin{lemma}}
	\newcommand{\elemma}{\end{lemma}}
\newcommand{\bnotation}{\begin{notation}}
	\newcommand{\enotation}{\end{notation}}
\newcommand{\bproof}{\begin{proof}}
	\newcommand{\eproof}{\end{proof}}
\newcommand{\bremark}{\begin{remark}}
	\newcommand{\eremark}{\end{remark}}
    \newcommand{\bremarks}{\begin{remarks}}
	\newcommand{\eremarks}{\end{remarks}}
\newcommand{\bcorollary}{\begin{corollary}}
	\newcommand{\ecorollary}{\end{corollary}}
\newcommand{\btheorem}{\begin{theorem}}
	\newcommand{\etheorem}{\end{theorem}}
 \newcommand{\bproblem}{\begin{problem}}
	\newcommand{\eproblem}{\end{problem}}
\newcommand{\bdefinition}{\begin{definition}}
	\newcommand{\edefinition}{\end{definition}}
\newcommand{\bequation}{\begin{equation}}
	\newcommand{\eequation}{\end{equation}}
\newcommand{\bequationn}{\begin{equation*}}
	\newcommand{\eequationn}{\end{equation*}}
\newcommand{\beqnarray}{\begin{eqnarray}}
	\newcommand{\eeqnarray}{\end{eqnarray}} \newcommand{\beqnarrayn}{\begin{eqnarray*}}
	\newcommand{\eeqnarrayn}{\end{eqnarray*}}
 \newcommand{\bexample}{\begin{example}}
	\newcommand{\eexample}{\end{example}}
     \newcommand{\bexamples}{\begin{examples}}
	\newcommand{\eexamples}{\end{examples}}
  \newcommand{\bassumption}{\begin{assumption}}
	\newcommand{\eassumption}{\end{assumption}}
\newcommand{\loc}{{\rm loc}}
\numberwithin{equation}{section}

\newcommand\blue[1]{\textcolor{blue}{#1}}

\newcommand{\Hmm}[1]{\leavevmode{\marginpar{\tiny%
			$\hbox to 0mm{\hspace*{-0.5mm}$\leftarrow$\hss}%
			\vcenter{\vrule depth 0.1mm height 0.1mm width \the\marginparwidth}%
			\hbox to
			0mm{\hss$\rightarrow$\hspace*{-0.5mm}}$\\\relax\raggedright #1}}}

\DeclareMathOperator{\diam}{diam}
\DeclareMathOperator{\supp}{supp}
\DeclareMathOperator{\capacity}{Cap}
\DeclareMathOperator{\vol}{Vol}
\DeclareMathOperator{\sgn}{sgn}
\newcommand{\core}{C_c^{\infty}(\Omega)}

\newcommand{\dHnn}{\,\mathrm{d}\mathcal{H}^{n-1}}
\newcommand{\dx}{\,\mathrm{d}x}

\newcommand{\dt}{\,\mathrm{d}t}

\newcommand{\ds}{\,\mathrm{d}s}
\newcommand{\dr}{\,\mathrm{d}r}

\newcommand{\dsigma}{\,\mathrm{d}\sigma}

\newcommand{\dvrho}{\,\mathrm{d}\varrho}

\DeclareMathOperator{\dive}{div}

\def\<{\langle}
\def\>{\rangle}

\long\def\prob#1\soln#2\endps{{\color{blue}#1}\medskip\par
	\noindent\underline{\sc Solution}:\hspace*{1em}\parindent=2em #2}

\AtBeginDocument{\begin{CJK}{UTF8}{gbsn}}
	\AtEndDocument{\end{CJK}}

            \def\gg{\gamma}

      \def\gw{\omega}
                
\def\Gg{\Gamma}

\def\Gw{\Omega}              
\pgfplotsset{compat=1.18}
\begin{document}
	\pagenumbering{gobble}
	\title{\textbf{Optimal Hardy-weights for the Finsler $p$-Dirichlet integral with a potential}}
\author{Yongjun Hou\thanks{Department of Mathematics, Technion -
Israel Institute of Technology, Haifa 3200003, Israel; yongjun.hou@campus.technion.ac.il; houmathlaw@outlook.com}}
 \date{December 24, 2025}
	\maketitle
	\pagenumbering{arabic}
	\vspace{-10mm}	
 	\begin{abstract} Fix an integer $n\geq 2$, an exponent $1<p<\infty$, and a domain $\Omega\subseteq\mathbb{R}^{n}$. Let $\Omega^{*}\triangleq\Omega\setminus\{\hat{x}\}$ where $\hat{x}\in\Omega$. Under some further conditions, we construct optimal Hardy-weights for the Finsler $p$-Dirichlet integral
$$Q_{0}[\phi;\Omega^{*}]\triangleq\int_{\Omega^{*}}H(x,\nabla \phi)^{p}\dx\quad \mbox{on}\quad C^{\infty}_{c}(\Omega^{*}),$$
and the Finsler $p$-Dirichlet integral with a potential
$$Q_{V}[\phi;\Omega]\triangleq\int_{\Omega}\left(H(x,\nabla \phi)^{p}+ V|\phi|^{p}\right)\dx\quad \mbox{on}\quad C^{\infty}_{c}(\Omega),$$where $H(x,\cdot)$ is a family of norms on $\mathbb{R}^{n}$ parameterized by $x\in\Omega^{*}$ or $x\in\Omega$, respectively, and the potential $V$ lies in a subspace $\widehat{M}^{q}_{\loc}(p;\Omega)$ of a local Morrey space $M^{q}_{\loc}(p;\Omega)$.
\end{abstract}
\medskip

 \section{Introduction} 
We start with the classical one-dimensional and multi-dimensional Hardy inequalities.
\btheorem [{\cite[Corollary 1.2.2]{Balinsky}}]
Let~$1<p<\infty$ and let~$f$ be a locally absolutely continuous function on~$(0,\infty)$ such that~$f'\in L^{p}((0,\infty))$ and that~$\lim_{x\rightarrow0^{+}}f(x)=0$. Then
$$\left(\frac{p-1}{p}\right)^{p}\int_{0}^{\infty}\frac{|f(x)|^{p}}{|x|^{p}}\dx\leq \int_{0}^{\infty}|f'(x)|^{p}\dx.$$
 The constant~$((p-1)/p)^{p}$ is sharp and the equality can only be attained by the zero function. 
\etheorem 
 \btheorem [{\cite[Corollary 1.2.6 and p.~12]{Balinsky}}]\label{multihardy}
 \begin{itemize}
\item For all functions
\begin{align*}
\phi\in \begin{cases} 		C^{\infty}_{c}(\R^{n}\setminus\{0\}) &~\mbox{if}~n<p<\infty,\\	
C^{\infty}_{c}(\R^{n})&~\mbox{if}~1\leq p<n,
 \end{cases}
 \end{align*}
 it holds that
 $$ \left|\frac{p-n}{p}\right|^{p}\int_{\R^{n}}\frac{|\phi(x)|^{p}}{|x|^{p}}\dx\leq\int_{\R^{n}}|\nabla \phi(x)|^{p}\dx,$$
 where the constant~$\left|(p-n)/p\right|^{p}$ is sharp.
 \item  For $p=n$, there is no positive constant~$C$ such that for all~$\phi\in C^{\infty}_{c}(\R^{n})$,  $$C\int_{\R^{n}}\frac{|\phi(x)|^{n}}{|x|^{n}}\dx\leq\int_{\R^{n}}|\nabla \phi(x)|^{n}\dx.$$
 \end{itemize}
  \etheorem
  \textbf{Throughout this paper, unless otherwise stated, we fix an integer~$n\geq 2$, an exponent~$1<p<\infty$, and a \emph{domain} (i.e., a nonempty, connected, open set) $\Omega\subseteq\R^{n}$, and denote~$\Omega\setminus\{\hat{x}\}$ by~$\Omega^{*}$ where~$\hat{x}\in\Omega$.}
  
  Hardy-type inequalities are a profound tool with important applications in mathematical analysis. For example, the monograph \cite{Kutev} by Kutev and Rangelov covers some applications to singular and degenerate heat and $p$-heat equations in its concluding chapter. Hardy-type inequalities are studied across multiple contexts and from diverse perspectives. For some specific works, see, e.g., \cite{Balinsky,Das2,Davies,Goel,PLMS,Kovarik,Lamberti,Shafrir} and references therein. 
\begin{definition}\label{efdfn}
 \emph{The \emph{energy functional}~$Q$ on~$\core$ is defined by \begin{align*}
Q[\phi]\triangleq Q[\phi;\Omega]\triangleq Q_{V}[\phi]&\triangleq Q_{V}[\phi;\Omega]\\
&\triangleq\int_{\Omega}\left(H(x,\nabla \phi)^{p}+ V|\phi|^{p}\right)\dx\\
&=\int_{\Omega}\left(\mathcal{A}(x,\nabla \phi)\cdot\nabla\phi+ V|\phi|^{p}\right)\dx,
\end{align*}
 where $H(x,\cdot)$, parameterized by~$x\in\Omega$, is a family of norms on $\mathbb{R}^{n}$ satisfying Assumptions~\ref{ass9}, the operator~$\mathcal{A}(x,\xi)\triangleq\nabla_{\xi}\left(H(x,\xi)^{p}/p\right)$ for almost all~$x\in\Omega$ and all~$\xi\in\mathbb{R}^{n}$, and the \emph{potential} $V$ belongs to a local Morrey space $M^{q}_{\loc}(p;\Omega)$ (see Definition~\ref{Morreydef1}).}
 \end{definition}
 \begin{definition}
\emph{For each nonempty open set~$D\subseteq\Omega$, the functional~$Q$ is called \emph{nonnegative} in~$D$ if for all~$\phi\in C^{\infty}_{c}(D)$, $Q[\phi;D]\geq 0$. When $V=0$, $Q_{0}$ is called the \emph{Finsler $p$-Dirichlet integral}.}
		\end{definition}
        Note that by \cite[Remark 2.13]{Hou2} (see also \cite[Lemma 4.18]{HPR}), if~$Q$ is nonnegative in~$\Omega$, then~$Q[\phi]\geq 0$ for all~$\phi\in W^{1,p}(\Omega)\cap C_{c}(\Omega)$.

        In the study of Hardy-type inequalities, Hardy-weights are an important concept. For example, in Theorem \ref{multihardy}, for~$n<p<\infty$ and~$1\leq p<n$,~$1/|x|^{p}$ is called a \emph{Hardy-weight} of the~\emph{$p$-Dirichlet integral}~$\int|\nabla\phi|^{p}\dx$ in~$\R^{n}\setminus\{0\}$ and~$\R^{n}$, respectively. In our setting, see the definition below.       \begin{definition}\label{hweights}
 \emph{A function~$g\in L^{1}_{\loc}(\Omega)$ is called a \emph{Hardy-weight of~$Q$} in $\Gw$ if there exists a constant~$C>0$ such that for all~$\phi\in\core$, the following \emph{Hardy-type inequality} holds:\bequationn\label{Hineq}
Q[\phi] \geq C\int_{\Omega}|g||\phi|^{p}\dx.
\eequationn
} 
\end{definition}
\blemma[{\cite[Lemma~2.14]{Hou2}}]
All Hardy-weights of~$Q$ in $\Gw$ form a linear space, which is denoted by~$\mathcal{H}(\Omega)$.
\elemma
The Hardy constant is defined as follows.
\bdefinition\label{hardyconstant}
\emph{Suppose that~$Q$ is nonnegative in~$\Omega$. For every~$g\in\mathcal{H}(\Omega)$, the \emph{(best) Hardy constant} $S_{g}$ is defined as $$\inf\bigg\{\frac{Q[\phi]}{\int_{\Omega}|g||\phi|^{p}\dx}~\bigg|~\phi\in W^{1,p}(\Omega)\cap C_{c}(\Omega)~\mbox{and} \int_{\Omega}|g||\phi|^{p}\dx>0\bigg\}.$$}
\edefinition
For a Maz'ya-type characterization of Hardy-weights, see \cite{Das,Hou2}, where the authors also established some attainments of the Hardy constant in a generalized Beppo Levi space or a subspace~$\widetilde{W}^{1,p}_{0}(\Omega)$ of a~$Q$-Sobolev space, respectively.

With regard to his celebrated estimates of solutions of second-order elliptic equations, Agmon \cite{Agmon} mentioned a research topic of looking for Hardy-weights that is as large as possible for an elliptic operator in some neighborhood of infinity in~$\R^{n}$. In recent years, relevant substantial achievements, called \emph{optimal Hardy-weights}, have been made for diverse operator classes. The pioneering work \cite{Devyver2} by Devyver, Fraas, and Pinchover dealt with second-order \emph{linear} elliptic operators. In the absence of a potential, optimal Hardy-weights were designed for the $p$-Laplacian and its weighted generalization, the~$(p,A)$-Laplacian, by Devyver and Pinchover \cite{Devyver1} (see Theorem~\ref{ohp} below) and Versano \cite{Versano}, respectively. In \cite{Versano}, under certain conditions, Versano also constucted optimal Hardy-weights for the~$(p,A)$-Laplacian with a potential. Discrete contexts have also witnessed some relevant advances. Keller, Pinchover, and Pogorzelski \cite{Keller1}, and Fischer \cite{Fishero} obtained optimal Hardy-weights for Schrödinger operators and for a quasilinear counterpart on graphs, respectively. We also mention the relevant works \cite{Dasoh,ohwfl} concering fractional Laplacians on the discrete half-line and on the integers, respectively, and \cite{Lemm} concerning the Euclidean lattice. For more discoveries on optimal Hardy-weights, see, e.g., \cite{Miura} by Miura, \cite{PVersano} by Pinchover and Versano, and \cite{Takeda} by Takeda.

It can be easily checked (cf.~\cite[Definition~6.3]{HPR}) that the following definition of criticality/subcriticality/supercriticality is equivalent to \cite[Definition~6.3]{HPR}. 
\begin{Def}
  	\ {\em 
  	   	  \begin{itemize}
    \item If~$\mathcal{H}(\Omega)=\{0\}$, 
the functional $Q$ is called \emph{critical} in $\Gw$;
    \item if~$\mathcal{H}(\Omega)\neq\emptyset$ and~$\mathcal{H}(\Omega)\neq\{0\}$, the functional $Q$ is called  \emph{subcritical} in $\Gw$;
    \item if~$\mathcal{H}(\Omega)=\emptyset$, the functional~$Q$ is called \emph{supercritical} in $\Gw$.
  \end{itemize}
}
\end{Def}
\begin{Def}
   \emph{Suppose that $Q$ is nonnegative in $\Omega$. A nonnegative sequence~$\{\phi_{k}\}_{k\in\mathbb{N}}\subseteq W^{1,p}(\Omega)\cap C_{c}(\Gw)$ is called a \emph{null-sequence} with respect to $Q$ in~$\Omega$ if
    \begin{itemize}
      \item there exists a fixed open set~$U\Subset\Omega$ such that~$\Vert \phi_{k}\Vert_{L^{p}(U)}=1$ for all~$k\in\mathbb{N}$;
      \item $\displaystyle{\lim_{k\rightarrow\infty}}Q[\phi_{k}]=0.$
    \end{itemize}}
  \end{Def}
  It turns out that under certain assumptions, a nonnegative functional~$Q$ is critical in~$\Omega$ if and only if~$Q$ admits a null-sequence (see Theorem \ref{nullc} (i)).
  \begin{Def}
  \emph{Suppose that $Q$ is nonnegative in $\Omega$. A \emph{ground state} of $Q$ is a positive function~$\phi\in W^{1,p}_{\loc}(\Omega)$ which is an $L^{p}_{\loc}(\Omega)$ limit of a null-sequence.}
  \end{Def}

In \cite{Devyver1}, optimal Hardy-weights for the~$p$-Dirichlet integral are defined as follows.
\bdefinition\label{ohwpd}
\emph{A nonnegative function~$W\in L^{\infty}_{\loc}(\Omega^{*})\setminus\{0\}$ is an \emph{optimal Hardy-weight} of the \emph{$p$-Dirichlet integral}~$\int_{\Omega^{*}}|\nabla\phi|^{p}\dx$ if
\begin{itemize}
\item the functional~$\int_{\Omega^{*}}|\nabla\phi|^{p}\dx-\int_{\Omega^{*}}W|\phi|^{p}\dx$ is critical in~$\Omega^{*}$; 
\item the functional~$\int_{\Omega^{*}}|\nabla\phi|^{p}\dx-\int_{\Omega^{*}}W|\phi|^{p}\dx$ is null-critical at~$\hat{x}$ and at infinity in the following sense: For any precompact open set~$O$ containing~$\hat{x}$ with~$\overline{O}\subseteq\Omega$, the ground state~$v$ of~$\int_{\Omega^{*}}|\nabla\phi|^{p}\dx-\int_{\Omega^{*}}W|\phi|^{p}\dx$ in~$\Omega^{*}$ satisfies:
$$\int_{O\setminus\{\hat{x}\}}|\nabla v|^{p}\dx=\infty\quad\mbox{and}\quad\int_{\Omega\setminus\overline{O}}|\nabla v|^{p}\dx=\infty;$$
\item $1$ is also the best constant~$\lambda$ for the inequality $$\int_{\Omega^{*}}|\nabla\phi|^{p}\dx\geq\lambda\int_{\Omega^{*}}W|\phi|^{p}\dx,$$ restricted to~$C^{\infty}$ functions~$\phi$ that are compactly supported either in a fixed punctured neighborhood of~$\hat{x}$ or in a fixed neighborhood (not containing~$\hat{x}$) of infinity in~$\Omega$.
\end{itemize}}
\edefinition
Suppose that for a nonnegative function~$W\in L^{\infty}_{\loc}(\Omega^{*})\setminus\{0\}$, $\int_{\Omega^{*}}|\nabla\phi|^{p}\dx-\int_{\Omega^{*}}W|\phi|^{p}\dx$ is critical in~$\Omega^{*}$ with a ground state~$v$. By \cite[Corollary 3.4]{Kovarik}, if~$\int_{\Omega^{*}}Wv^{p}\dx=\infty$, then $1$ is also the best constant~$\lambda$ for the inequality $$\int_{\Omega^{*}}|\nabla\phi|^{p}\dx\geq\lambda\int_{\Omega^{*}}W|\phi|^{p}\dx,$$ restricted to~$C^{\infty}$ functions~$\phi$ that are compactly supported in any fixed neighborhood of infinity in~$\Omega^{*}$. See also Lemma \ref{optimalinf}.

The following theorem is the constructed optimal Hardy-weights for the $p$-Dirichlet integral $\int_{\Omega^{*}}|\nabla\phi|^{p}\dx$ in \cite{Devyver1}. We define~$f_{0}(t)\triangleq t^{(p-1)/p}$ for all~$t\in(0,\infty)$.
\btheorem[{\cite[Theorem~1.5]{Devyver1}}]\label{ohp}
Suppose that there is a positive~$p$-harmonic function~$\mathcal{G}$ in~$\Omega^{*}$ satisfying, for some nonnegative real number~$\sigma$,
 \begin{align*}
\begin{cases}		\lim_{x\rightarrow \hat{x}}\mathcal{G}(x)=\infty\quad\mbox{and}\quad\lim_{x\rightarrow \overline{\infty}}\mathcal{G}(x)=0&\mbox{if}~1<p\leq n,\\
\lim_{x\rightarrow \hat{x}}\mathcal{G}(x)=\sigma\quad\mbox{and}\quad\lim_{x\rightarrow \overline{\infty}}\mathcal{G}(x)=\begin{cases}
    \infty&~ \mbox{if}~\sigma=0,\\    0&~ \mbox{if}~\sigma>0,
\end{cases}&\mbox{if}~p>n.
		\end{cases}
\end{align*}
Let~$W_{0}\triangleq\left(\frac{p-1}{p}\right)^{p}\left|\frac{\nabla\mathcal{G}}{\mathcal{G}}\right|^{p}$ and
\begin{align*}
\begin{cases} 
v\triangleq f_{0}(\mathcal{G})~\mbox{and}~W\triangleq W_{0}&\mbox{if}~1<p\leq n, \mbox{or}~p>n~\mbox{and}~\sigma=0,\\ v\triangleq f_{0}(\mathcal{G}(\sigma-\mathcal{G}))~\mbox{and}~W\triangleq\\\left(\frac{1}{\sigma-\mathcal{G}}\right)^{p}W_{0}|\sigma-2\mathcal{G}|^{p-2}(2(p-2)\mathcal{G}(\sigma-\mathcal{G})+\sigma^{2})&\mbox{otherwise}.
\end{cases}
\end{align*}
Then~$W$ is an optimal Hardy-weight of the $p$-Dirichlet integral~$\int_{\Omega^{*}}|\nabla\phi|^{p}\dx$ in the sense of Definition \ref{ohwpd}.
\etheorem
\begin{definition}\label{BREGMAN1}
	\emph{Let~$f: X\rightarrow \mathbb{R}$ be a G\^{a}teaux differentiable convex function on a Banach space~$(X,\Vert\cdot\Vert)$, and denote by $T_x$ its G\^{a}teaux derivative at $x\in X$. Then the function
		\begin{align*}
        D_{f}:X\times X&\longrightarrow\R,\\
		(y,x)&\longmapsto f(y)-f(x)-T_{x}(y-x),
		\end{align*}
		is called the \emph{Bregman distance} of~$f$.}
\end{definition}
In the constructions of \cite{Devyver1} and \cite{Versano}, the upper bounds for the corresponding Bregman distances in the following estimates are pivotal, while the lower bounds are related to criticality theory (see \cite{PPAPDE}).
\blemma[{\cite[(3.11)]{Regev}}]\label{ABE}
Let~$A\in \R^{n\times n}$ be a symmetric positive definite matrix. Then for all~$\xi,\eta\in\R^{n}$,
$$D_{|\cdot|_{A}^{p}}(\xi+\eta,\xi)=|\xi+\eta|_{A}^{p}-|\xi|_{A}^{p}-p|\xi|^{p-2}A\xi\cdot\eta\asymp |\eta|_{A}^{2}(|\eta|_{A}+|\xi|_{A})^{p-2},$$
where~$|\xi|_{A}\triangleq\sqrt{A\xi\cdot\xi}$ and both equivalence constants depend only on~$p$.
\elemma
In the present paper, building on the works \cite{Devyver1} and \cite{Versano}, under certain hypotheses, we produce optimal Hardy-weights (see Definition~\ref{ohw}) for the Finsler $p$-Dirichlet integral
$$Q_{0}[\phi;\Omega^{*}]\triangleq\int_{\Omega^{*}}H(x,\nabla \phi)^{p}\dx\quad \mbox{defined on}\quad C^{\infty}_{c}(\Omega^{*}),$$
and the Finsler $p$-Dirichlet integral with a potential
$$Q_{V}[\phi;\Omega]\triangleq\int_{\Omega}\left(H(x,\nabla \phi)^{p}+ V|\phi|^{p}\right)\dx\quad \mbox{defined on}\quad C^{\infty}_{c}(\Omega),$$ where the norm family~$H(x,\cdot)$ satisfies some assumptions and the potential $V$ belongs to a certain local Morrey space. The local Euler-Lagrange equation of~$Q_{V}[\cdot;\Omega]$ in~$\Omega$ is the quasilinear elliptic partial differential equation
 $$Q'[u]\triangleq-\mathrm{div}\mathcal{A}(x,\nabla u)+ V|u|^{p-2}u=0\quad\mbox{in}~\Omega.$$For~$V=0$, the equation~$Q'[u]=0$ reduces to the \emph{Finsler~$p$-Laplace equation} (see, e.g., \cite{Combete}), which, in~$\Omega^{*}$, is associated with the former functional~$Q_{0}[\cdot;\Omega^{*}]$. See Sections~\ref{avint} and~\ref{fqeq} for more details. For more insights into the Finsler~$p$-Laplace equation, see, e.g., \cite{Bianchini,Jaros,Jarosalap}. 
  
  As in \cite{Devyver1} and \cite{Versano}, the Bregman distances of~$H(x,\cdot)^{p}$ are also involved in our setting. In both cases, as a key ingredient, we \emph{assume} an upper bound (see Assumption~\ref{assup}) for the Bregman distances of~$H(x,\cdot)^{p}$, uniform with respect to almost all~$x\in\Omega^{*}$ or almost all~$x\in\Omega$. We also present a class of norms (see Example~\ref{examplesu}) realizing this assumption and the others. For more on Bregman distances, we refer the readers to, e.g., \cite{Information,Lindqvist,Pinliou,Shafrir,Sprung}.
    
     In the case of zero potential, we use positive Finsler $p$-harmonic functions (see Definition~\ref{aharmonic}) satisfying some boundary conditions, while for nonzero potentials, we exploit Green potentials (see Definition~\ref{gpotential}) with some extra properties. Moreover, we provide some specific examples of Finsler $p$-harmonic functions and Green potentials, with the corresponding properties. In some examples, the principal eigenvalue of~$Q'$ is involved. Thus we also study eigenvalues of~$Q'$ in domains compactly contained in~$\Omega$. Under only Assumptions~\ref{ass9} and for~$V\in M^{q}_{\loc}(p;\Omega)$, we demonstrate that the principal eigenvalue of~$Q'$ is isolated and that all the real eigenvalues of~$Q'$ form a closed subset of~$\R$.

     For clarity, we present a special case of the derived optimal Hardy-weights without a potential, which can be deduced directly from Theorem \ref{zerop} and Examples \ref{fexas} (iii). 
     \bexample
Suppose that~$p\neq n$. Let~$\Omega=\R^{n}$ and let~$H$ be a fixed norm on~$\R^{n}$ satisfying certain assumptions. For all~$x\in\R^{n}$, let$$H_{0}(x)\triangleq \sup_{\xi\in\R^{n}\setminus\{0\}}\frac{x\cdot\xi}{H(\xi)}$$ be the dual norm of~$H$ and let~$\mathcal{G}(x)\triangleq H_{0}(x)^{(p-n)/(p-1)}$ for~$x\neq 0$. Then in~$\R^{n}\setminus\{0\}$,~$W\triangleq\left(\frac{p-1}{p}\right)^{p}H\left(\frac{\nabla\mathcal{G}}{\mathcal{G}}\right)^{p}$ is an optimal Hardy-weight (in the sense of Definition \ref{ohw}) of~$Q_{0}$ and (up to a positive multiplicative constant)~$H_{0}^{(p-n)/p}$ is the unique ground state of~$Q_{-W}$.
     \eexample

The paper is organized as follows. In Section~\ref{equationf}, we introduce the norm family~$H$, the operator~$\mathcal{A}$, and the local Morrey space $M^{q}_{\loc}(p;\Omega)$ together with its subspace~$\widehat{M}^{q}_{\loc}(p;\Omega)$. We also discuss the equation~$Q'[u]=0$ and review some results and a notion of criticality theory. In Section \ref{zpotential}, we define and investigate optimal Hardy-weights. In Sections~\ref{zeropcons} and~\ref{nzpotential}, we construct optimal Hardy-weights for the functional~$Q$ with zero and nonzero potentials, respectively. In the appendix, we establish the two results on eigenvalues of~$Q'$.
\section*{Basic symbols}
 \normalsize
 \vspace{-1mm}
\begin{longtable}[1]{p{60pt} p{350pt} }
 $\mathbb{N}$ & the set of positive integers\\
 $c_{p}$& the positive constant~$(p/(p-1))^{p-1}$\\
 $\vol(\Gg)$ & the Lebesgue measure of a measurable subset~$\Gg$ of~$\R^{n}$\\
  $\mathcal{H}^{n-1}$ & the ($n-1$)-dimensional Hausdorff measure on~$\R^{n}$\\
     $\sgn(\cdot)$& the sign function on~$\R$\\ 
    $|\xi|$& the Euclidean norm of a vector~$\xi\in\R^{n}$\\
    $|\xi|_s$& the $s$-norm of a vector~$\xi=(\xi_{1},\xi_{2},\ldots,\xi_{n})\in\R^{n}$ ($1<s<\infty$), i.e.,~$|\xi|_{s}\triangleq\left(\sum_{i=1}^{n}|\xi_{i}|^{s}\right)^{1/s}$\\
    $s'$ & the conjugate exponent~$s/(s-1)$ of~$s\in(1,\infty)$\\
 $\overline{\infty}$ & unless otherwise stated, the ideal point in the one-point compactification of~$\Omega$\\
 $\mathring{K},\overline{K}$ & the interior and closure of~$K\subseteq\R^{n}$, respectively \\
 $\mathrm{diam}(K)$ & the diameter of~$K\subseteq\R^{n}$\\
	$E \Subset F$ &$\overline{E}$ is a compact subset of $F$ with~$E,F\subseteq\R^{n}$ \\
	$B_r(x)$  & the open ball centered at $x\in\mathbb{R}^{n}$ with radius $r>0$\\
    $S_r(x)$  & the sphere centered at $x\in\mathbb{R}^{n}$ with radius $r>0$\\
$\mathcal{F}_{c}(U)$ & the space~$\{f\in\mathcal{F}(U)~|~\supp f\Subset U\}$, where~$U\subseteq\R^{n}$ is a nonempty open set and~$\mathcal{F}(U)$ is an arbitrary linear space of functions on~$U$\\
 $f^{+}$  & the positive part of a function $f$, namely $\max\{f,0\}$\\
 $f^{-}$  & the negative part of a function $f$, namely $-\min\{f,0\}$\\
$C$ & a positive constant which may vary from line to line\\
$f\asymp g$&$cg\leq f\leq Cg$ for some positive constants~$c$ and~$C$ (called \emph{equivalence constants}), where $f$ and $g$ are nonnegative functions and in this case, are called \emph{equivalent}\\
		\end{longtable}
\section{Finsler~$p$-Laplace equation with a potential}\label{equationf}
In this section, after specifying the norm family~$H$ and the potential~$V$, we explicitly present the Finsler $p$-Laplace equation with a potential~$Q'[u]=0$. We also  review some criticality theory, laying a foundation for Sections~\ref{zpotential}, \ref{zeropcons}, and~\ref{nzpotential}.
	\subsection{The norm family~$H$, the operator~$\mathcal{A}$, and the potential~$V$}\label{avint}
In this subsection, we first define the norm family~$(H(x,\cdot))_{x\in\Omega}$ on $\R^{n}$ and subsequently derive from~$H$ an associated variational Lagrangian~$F$ and the operator~$\mathcal{A}$ . Next, we review the local Morrey space~$M^{q}_{\loc}(p;\Omega)$ and its subspace~$\widehat{M}^{q}_{\loc}(p;\Omega)$. Furthermore, we impose an additional assumption on~$H$, which, combined with the condition~$V\in\widehat{M}^{q}_{\loc}(p;\Omega)$, ensures the continuous differentiability of solutions to some relevant equations.

Our discussion starts with the following norm family~$H$. See also \cite[Assumptions~2.3]{Hou2}.	\begin{assumptions}\label{ass9} 
		{\em			Let~$H:\Omega\times\R^{n}$ be a mapping such that for almost all~$x\in\Omega$,~$H(x,\cdot)$ is a norm on $\mathbb{R}^{n}$, with the following properties:   
			\begin{itemize}
				\item {\bf (Measurability)} For all~$\xi\in\mathbb{R}^{n}$, the mapping $x\mapsto H(x,\xi)$ is measurable in $\Gw$;
				\item {\bf (Local uniform equivalence)} for every domain $\omega\Subset \Gw$, there exist two constants~$0<\kappa_\omega\leq\nu_\omega<\infty$ such that for almost all $x\in \omega$ and all $\xi\in \R^n$,$$\kappa_\omega |\xi|\leq H(x,\xi)\leq\nu_\omega|\xi|;$$
				\item {\bf (Uniform convexity)} for almost all~$x\in\Omega$, the Banach space $\mathbb{R}^{n}_{x}\triangleq(\mathbb{R}^{n},H(x,\cdot))$ is uniformly convex; 
				\item {\bf (Differentiability with respect to $\xi$)} for almost all~$x\in \Gw$, the mapping $\xi\mapsto H(x,\xi)$ is differentiable in $\R^n\setminus\{0\}$.
			\end{itemize}		
		}
	\end{assumptions}
The norm family~$H$ induces the Lagrangian~$F$. 
\begin{theorem}[{\cite[Theorem~2.4]{Hou2}}]
  Let~$(H(x,\cdot))_{x\in\Omega}$ be a family of norms on $\mathbb{R}^{n}$ fulfilling Assumptions~\ref{ass9} and for almost all
 ~$x\in \Omega$ and all~$\xi\in \R^{n}$, let
	$$F(x,\xi)\triangleq\frac{1}{p}H(x,\xi)^{p}.$$
            Then
			\begin{itemize}
				\item {\bf (Measurability)} for all~$\xi\in\mathbb{R}^{n}$, the mapping $x\mapsto F(x,\xi)$ is measurable in $\Gw$;
				\item {\bf (Local uniform ellipticity and boundedness)} 
				for every domain $\omega\Subset \Gw$, there exist two constants $0<\kappa'_\omega\leq\nu'_\omega<\infty$ such that for almost all $x\in \omega$ and all $\xi\in \R^n$,$$\kappa'_\omega|\xi|^{p}\leq F(x,\xi)\leq\nu'_\omega|\xi|^{p};$$
				\item {\bf (Strict convexity and $C^1$ with respect to $\xi$)}  for almost all~$x\in \Gw$, the mapping $\xi\mapsto F(x,\xi)$ is strictly convex and continuously differentiable in $\R^n$;
				\item {\bf (Homogeneity)} for almost all~$x\in \Gw$, all~$\lambda\in\mathbb{R}$, and all~$\xi\in\mathbb{R}^{n}$, $F(x,\lambda\xi)=|\lambda|^{p}F(x,\xi)$.
			\end{itemize}		
	\end{theorem}
    Note that in \cite{HPR}, the Lagrangian~$F$ was adopted as the basis of that work (see \cite[Assumption~2.1]{HPR}). However, the two ways are equivalent (see \cite[Remark~2.5]{Hou2}). 
    
Next we define and discuss the operator~$\mathcal{A}$.	\begin{definition}
\emph{Let~$(H(x,\cdot))_{x\in\Omega}$ be a family of norms on $\mathbb{R}^{n}$ fulfilling Assumptions~\ref{ass9}. For almost all~$x\in\Omega$ and all~$\xi\in\mathbb{R}^{n}$, let
			$$ \mathcal{A}(x,\xi) \triangleq \nabla_\xi F(x,\xi).$$
            }
	\end{definition}
    \begin{notation}
\emph{Let~$(H(x,\cdot))_{x\in\Omega}$ be a family of norms on $\mathbb{R}^{n}$ fulfilling Assumptions~\ref{ass9}. For almost all~$x\in\Omega$ and all~$\xi\in\mathbb{R}^{n}$, we also denote~$H(x,\xi)$ by~$|\xi|_{\mathcal{A}}$ in some places. When the norm family~$(H(x,\cdot))_{x\in\Omega}$ is independent of~$x$, with a slight abuse of notation, we simply write~$H(x,\xi)$ and $\mathcal{A}(x,\xi) $  as~$H(\xi)$ and~$\mathcal{A}(\xi)$, respectively.}
	\end{notation}
The following relationship between~$F$ and $\mathcal{A}$ is due to the homogeneity of $F$.
	\begin{lemma}[{\cite[p.~100]{HKM}}]\label{pflemma}	Let~$(H(x,\cdot))_{x\in\Omega}$ be a family of norms on $\mathbb{R}^{n}$ fulfilling Assumptions~\ref{ass9}. For almost all~$x\in\Omega$ and all~$\xi\in\mathbb{R}^{n}$,~$\mathcal{A}(x,\xi)\cdot\xi =pF(x,\xi)=H(x,\xi)^{p}.$
	\end{lemma}
 The following theorem consists of some fundamental properties of the operator~$\mathcal{A}$.
		\begin{Thm}[{\cite[Lemma~5.9]{HKM}}]\label{thm_1}
			Let~$(H(x,\cdot))_{x\in\Omega}$ be a family of norms satisfying Assumptions~\ref{ass9}. For every domain $\omega\Subset\Omega$, let  $\alpha_{\omega}=\kappa'_{\omega}$, $\beta_{\omega}=2^{p}\nu'_{\omega}$. Then the vector-valued function~$\mathcal{A}:  \Gw\times \mathbb{R}^{n}\rightarrow \mathbb{R}^{n}$ satisfies the following properties:
				\begin{itemize}
					\item {\bf (Continuity and measurability)} For almost all $x\in \Gw$, the function $\xi\mapsto \mathcal{A}(x,\xi )$ is continuous in~$\R^{n}$, and for all~$\xi\in \mathbb{R}^{n}$, the function $x \mapsto \mathcal{A}(x,\xi)$ is Lebesgue measurable on $\Gw$;
					
					%
					\item {\bf (Local uniform ellipticity and boundedness)} for all domains $\omega\Subset \Gw$, almost all $x\in \omega$, and all $\xi \in \mathbb{R}^{n}$,		\begin{equation*}\label{structure}
						\alpha_\omega|\xi|^{p}\le\mathcal{A}(x,\xi)\cdot\xi\quad\mbox{and}\quad
						|\mathcal{A}(x,\xi)|\le \beta_\omega\,|\xi|^{{p}-1};
					\end{equation*}
					\item {\bf (Strict monotonicity)} for almost all~$x\in\Gw$ and all~$\xi,\eta\in\mathbb{R}^{n}$ ($\xi\neq\eta$),
					$$\big(\mathcal{A}(x,\xi)-\mathcal{A}(x,\eta)\big) \! \cdot \! (\xi-\eta)> 0;$$
                    \item {\bf (Homogeneity)} for all~$\lambda\in {\mathbb{R}\setminus\{0\}}$,
					$\mathcal{A}(x,\lambda \xi)=\lambda\,|\lambda|^{p-2}\,\mathcal{A}(x,\xi)$.
			\end{itemize}
		\end{Thm}
Now we are in a position to present the local Morrey space~$M^{q}_{\loc}(p;\Omega)$ and its subspace~$\widehat{M}^{q}_{\loc}(p;\Omega)$. See also \cite[Definitions 2.9 and 2.11]{Hou2}.
  \begin{Def}\label{Morreydef1}
  \emph{Suppose that~$\omega\Subset\Omega$ is a domain and that $f$ is a real-valued measurable function on~$\omega$. 
				\begin{itemize}
					\item For $p<n$ and $q>n/p$, we define$$\Vert f\Vert_{M^{q}(p;\omega)}\triangleq \sup_{\substack{y\in\gw\\0<r<\diam(\gw)}}
					\frac{1}{r^{n/q'}}\int_{\omega\cap B_{r}(y)}|f|\dx$$ and $$M^{q}(p;\omega)\triangleq\{f\in L^{1}_{\loc}(\omega)~|~\Vert f\Vert_{M^{q}(p;\omega)}<\infty\};$$ 
					\item for $p=n$ and $q>n$ , we define$$\Vert f\Vert_{M^{q}(n;\omega)}\triangleq \sup_{\substack{y\in\gw\\0<r<\diam(\gw)}} \varphi_{q}(r)\int_{\omega\cap B_{r}(y)}|f|\dx,$$
					where $\varphi_{q}(r)\triangleq \left(\log\big(\mathrm{diam}(\omega)/r\big)\right)^{q/n'}$, and
$$M^{q}(n;\omega)\triangleq\{f\in L^{1}_{\loc}(\omega)~|~\Vert f\Vert_{M^{q}(n;\omega)}<\infty\};$$
					\item for $p>n$ and $q=1$, we define~$M^{q}(p;\omega)\triangleq L^{1}(\omega)$.
				\end{itemize}
    Finally, we define the \emph{local Morrey space} by $$M^{q}_{\loc}(p;\Omega)\triangleq \bigcap_{\substack{\omega\Subset\Omega\\\omega~\mbox{is a domain}}}M^{q}(p;\omega).$$
			}
		\end{Def}
\begin{definition}
\emph{We denote by~$\widehat{M}^{q}_{\loc}(p;\Omega)$ the space of all the functions $f$ in~$M^{q}_{\loc}(p;\Omega)$ satisfying the following conditions: When~$p<n$,~$q>n$; 
when~$p\geq n$, for some~$\theta\in (p-1,p)$ and all domains~$\omega\Subset\Omega$,
$$\sup_{\substack{y\in\omega\\0<r<\diam(\omega)}}\frac{1}{r^{n-p+\theta}}\int_{B_{r}(y)\cap\omega}|f|\dx<\infty.$$
}
\end{definition}
\bassumption\label{C1alpha}
\emph{For all~$x\in\Omega$,~$\mathcal{A}(x,\cdot)\in C^{1}(\R^{n}\setminus\{0\})$. For every domain~$\omega\Subset\Omega$, 
there exist positive constants~$\vartheta_{\omega}\leq 1, \Lambda_{1,\omega},\Lambda_{2,\omega}$, and~$\Lambda_{3,\omega}$ such that for all~$x,y\in\omega$, all~$\xi\in\R^{n}\setminus\{0\}$, and all~$\eta\in\R^{n}$,
\begin{align*}
D_{\xi}\mathcal{A}(x,\xi)\eta\cdot\eta\geq \Lambda_{1,\omega}|\xi|^{p-2}|\eta|^{2},\quad |D_{\xi}\mathcal{A}(x,\xi)|\leq \Lambda_{2,\omega}|\xi|^{p-2},
\end{align*}
and \begin{align*}
|\mathcal{A}(x,\xi)-\mathcal{A}(y,\xi)|\leq \Lambda_{3,\omega}|\xi|^{p-1}|x-y|^{\vartheta_{\omega}}.
\end{align*}
}
\eassumption
See Example \ref{examplesu} for a norm family realizing Assumptions \ref{ass9} and \ref{C1alpha}.

\textbf{Throughout this paper, unless otherwise stated (especially in the appendix), we assume that $(H(x,\cdot))_{x\in\Gw}$ is a family of norms on $\mathbb{R}^{n}$ satisfying Assumptions~\ref{ass9} and~\ref{C1alpha}, and that~$V\in \widehat{M}^{q}_{\loc}(p;\Omega)$}.

Assumption~\ref{C1alpha} implies the following lower bound (cf. \cite[Assumption~2.8]{HPR}) for the Bregman distance of~$|\cdot|_{\mathcal{A}}^{p}$. The proof is based on Taylor's formula with integral remainder. See also \cite[p.~163]{Lindqvist} for a similar calculation.
\begin{lemma}\label{bllem}
	{\em
 For every domain $\gw\Subset \Gw$, there exists  
		a positive constant $C_{\omega}$ (independent of~$x\in\omega$, or~$\xi,\eta\in\R^{n}$) such that for all~$x\in \gw$ and~$\xi,\eta\in \R^n$,  
	\begin{align*}	|\xi+\eta|^{p}_{\mathcal{A}}-|\xi|^{p}_{\mathcal{A}}-p\mathcal{A}(x,\xi)\cdot\eta\geq \begin{cases}
				C_{\omega}|\eta|^{p}&\mbox{if $p\geq 2$,}\\
				C_{\omega}|\eta|^{2}(|\xi|+|\eta|)^{p-2}&\mbox{if $p<2$.}
		\end{cases}
   \end{align*}	
	}
\end{lemma}
\subsection{The equation~$Q'[u]=0$ and a one-sided simplified energy}\label{fqeq}
In this subsection, we define the equation~$Q'[u]=0$ and introduce another assumption on~$H$ which leads to a one-sided simplified energy. The one-sided simplified energy is instrumental in deriving optimal Hardy-weights.

        We first consider the \emph{nonhomogeneous Finsler~$p$-Laplace equation with a potential}.
 \begin{Def}
			\emph{  
Let~$g\in M^{q}_{\loc}(p;\Omega)$. A function~$v\in W^{1,p}_{\loc}(\Omega)$ is a  (\emph{weak}) \emph{solution} of the equation
				\begin{equation*}\label{half}
					Q'[u]\triangleq Q'_{V}[u]\triangleq -\dive\mathcal{A}(x,\nabla u)+V|u|^{p-2}u=g,
				\end{equation*}
				in~$\Omega$ if for all~$\phi \in C_{c}^{\infty}(\Omega)$,$$\int_{\Omega}\mathcal{A}(x,\nabla v)\cdot \nabla \phi\dx+\int_{\Omega}V|v|^{p-2}v \phi\dx=\int_{\Omega}g\phi\dx.$$
    A function~$v\in W^{1,p}_{\loc}(\Omega)$ is a (\emph{weak}) \emph{supersolution} of the above equation
				in~$\Omega$ if for all nonnegative~$\phi \in C_{c}^{\infty}(\Omega)$, $$\int_{\Omega}\mathcal{A}(x,\nabla v)\cdot \nabla \phi\dx+\int_{\Omega}V|v|^{p-2}v \phi\dx\geq \int_{\Omega}g\phi\dx.$$ A function $v\in W^{1,p}_{\loc}(\Omega)$ is a (\emph{weak}) \emph{subsolution} of the above equation
				in~$\Omega$ if for all nonnegative~$\phi \in C_{c}^{\infty}(\Omega)$, $$\int_{\Omega}\mathcal{A}(x,\nabla v)\cdot \nabla \phi\dx+\int_{\Omega}V|v|^{p-2}v \phi\dx\leq \int_{\Omega}g\phi\dx.$$
            A solution/supersolution/subsolution of the homogeneous equation~$Q'[u]=0$ in $\Omega$ is called \emph{$Q$-harmonic/$Q$-superharmonic/$Q$-subharmonic} in $\Omega$, respectively.}		\end{Def}
                \bdefinition\label{danames}
                \emph{The operator $\dive\mathcal{A}(x,\nabla u)$ is called the \emph{Finsler~$p$-Laplacian} or the \emph{$\mathcal{A}$-Laplacian}. If the norm family~$H$ is a fixed norm, $\dive\mathcal{A}(\nabla u)$ is also called the \emph{anisotropic $p$-Laplacian}.}
                \edefinition
    \bremark\label{condiffe}
\emph{Recall that for all nonnegative~$\phi\in\core$, solutions of the equation
\begin{align}\label{neqn}
-\dive\mathcal{A}(x,\nabla u)+V|u|^{p-2}u=\phi
\end{align}
are locally H\"older continuous and hence locally bounded in~$\Omega$. See \cite[Theorem~4.11]{Maly97} and \cite[Theorem~7.5.7]{Serrin}.
By \cite[Theorem~5.3]{Lieberman93}, 
under Assumption~\ref{C1alpha} and for~$V\in\widehat{M}^{q}_{\loc}(p;\Omega)$, 
solutions of Equation \eqref{neqn}
are continuously differentiable in~$\Omega$.}
\eremark
The following assumed upper bound for the Bregman distance of~$|\cdot|^{p}_{\mathcal{A}}$ will play an important role in constructing optimal Hardy-weights. See also \cite[Assumption~8.23]{Hou}.
\begin{assumption}\label{assup}
	{\em 
There exists  
		a positive constant $C$ (independent of~$x\in\Omega$ and~$\xi,\eta\in\R^{n}$) such that for almost all~$x\in \Gw$ and all~$\xi,\eta\in \R^n$, 		$$|\xi+\eta|^{p}_{\mathcal{A}}-|\xi|^{p}_{\mathcal{A}}-p\mathcal{A}(x,\xi)\cdot\eta\leq C|\eta|^{2}_\mathcal{A}(|\xi|_{\mathcal{A}}+|\eta|_{\mathcal{A}})^{p-2}.$$
	}
\end{assumption}
Now we present a class of examples realizing Assumptions~\ref{ass9},~\ref{C1alpha}, and~\ref{assup} simultaneously.
\bexample\label{examplesu}
\emph{Let $2\leq s<\infty$ and recall that $|\xi|_{s}\triangleq\left(\sum_{i=1}^{n}|\xi_{i}|^{s}\right)^{1/s}$ for all~$\xi\in\R^{n}$. Then (see \cite[Corollary~3.4]{Hou2})~$|\cdot|_{s}^{p}\in C^{1}(\R^{n})$ and for all~$\xi\in\R^{n}$ and~$i=1,2,\ldots,n$,
$$\frac{\partial|\xi|_{s}^{p}}{\partial\xi_{i}}=p|\xi|_{s}^{p-s}|\xi_{i}|^{s-2}\xi_{i}.$$
By \cite[Lemma~3.12]{Hou2}, there exists a positive constant $C(p,s)$ such that for all~$\xi,\eta\in \R^n$, 		\begin{align}\label{sbreg}
        |\xi+\eta|^{p}_{s}-|\xi|^{p}_{s}-p|\xi|_{s}^{p-s}\sum_{i=1}^{n}|\xi_{i}|^{s-2}\xi_{i}\eta_{i}\leq C(p,s)|\eta|^{2}_s(|\xi|_{s}+|\eta|_{s})^{p-2}.
        \end{align} 
Let~$A\in\R^{n\times n}$ be a symmetric positive definite matrix. Recall that~$|\xi|_{A}=\sqrt{A\xi\cdot\xi}$ for all~$\xi\in\R^{n}$. 
        By \cite[(3.11)]{Regev}, for all~$\xi,\eta\in \R^n$, 	\begin{align}\label{abreg}	    
		|\xi+\eta|^{p}_{A}-|\xi|^{p}_{A}-p|\xi|_{A}^{p-2}A\xi\cdot\eta\asymp|\eta|^{2}_A(|\xi|_{A}+|\eta|_{A})^{p-2},
        \end{align}
        where both equivalence constants depend only on~$p$. Now let~$H(x,\xi)\triangleq H(\xi)\triangleq\sqrt[p]{|\xi|_{s}^{p}+|\xi|_{A}^{p}}$ ($2\leq s<\infty$) for all~$x\in\Omega$ and~$\xi\in\R^{n}$. Then~$H$ satisfies Assumptions~\ref{ass9} and~\ref{C1alpha} (cf. \cite[Example~3.5]{Hou2}). There exists a positive constant~$C_{0}=C_{0}(s,A)$ such that for all~$\xi\in\R^{n}$ (and all~$x\in\Omega$),
        $$|\xi|_{s}\leq H(\xi)\leq C_{0}|\xi|_{s}\quad\mbox{and}\quad|\xi|_{A}\leq H(\xi)\leq C_{0}|\xi|_{A}.$$
Then by \eqref{sbreg} and \eqref{abreg}, for all~$\xi,\eta\in \R^n$ (and all~$x\in \Gw$),
\begin{align*}
&|\xi+\eta|^{p}_{\mathcal{A}}-|\xi|^{p}_{\mathcal{A}}-p\mathcal{A}(\xi)\cdot\eta\\
&= |\xi+\eta|^{p}_{s}-|\xi|^{p}_{s}-p|\xi|_{s}^{p-s}\sum_{i=1}^{n}|\xi_{i}|^{s-2}\xi_{i}\eta_{i}+|\xi+\eta|^{p}_{A}-|\xi|^{p}_{A}-p|\xi|_{A}^{p-2}A\xi\cdot\eta\\
&\leq C(p,s)|\eta|^{2}_s(|\xi|_{s}+|\eta|_{s})^{p-2}+C(p)|\eta|^{2}_A(|\xi|_{A}+|\eta|_{A})^{p-2}\\
&\leq C|\eta|^{2}_\mathcal{A}(|\xi|_{\mathcal{A}}+|\eta|_{\mathcal{A}})^{p-2},
\end{align*}
where~$C=\begin{cases}		C(p,s)&~\mbox{if}~p\geq 2,\\
		C(p,s,A)&~\mbox{if}~p<2.
		\end{cases}$ Hence~$H$ also fulfills Assumption~\ref{assup}.
        }
\eexample
The following corollary is called the \emph{one-sided simplified energy}. See also \cite[Lemma 3.4]{Regev} and \cite[Corollary 2.28~(i)]{Versano}.
\begin{corollary}[{\cite[Lemma~8.28]{Hou}}]\label{seup}
Suppose that Assumption~\ref{assup} holds. Take any positive
 subsolution $v\in W^{1,p}_{\loc}(\Omega)$ of $Q'[f]=0$ and any nonnegative function~$u\in W^{1,p}_{\loc}(\Omega)$. If~$u^{p}/v^{p-1} \in W^{1,p}_{c}(\Omega)$, the product and chain rules for~$u^{p}/v^{p-1}$ holds, and~$vw$ satisfies the product rule for~$w\triangleq u/v$,
then
  $$Q[vw]\leq C\int_{\Omega}v^{2}|\nabla w|^{2}_{\mathcal{A}}\left(w|\nabla v|_{\mathcal{A}}+v|\nabla w|_{\mathcal{A}}\right)^{p-2}\dx,$$
  where~$C$ is as in Assumption~\ref{assup}.
\end{corollary}
\bremark
\emph{For some two-sided simplified energies, see \cite[Lemma~8.28]{Hou}, \cite[Lemmas 3.18, 3.19, and 3.20]{Hou2} and \cite[Lemma 2.2]{Pinliou}.}
\eremark
By virtue of Corollary~\ref{seup} and H\"older's inequality, we may easily obtain the following upper bound. See also \cite[Corollary~2.28~(ii)]{Versano}.
\begin{corollary}\label{co1}
Suppose that Assumption~\ref{assup} holds. Let~$v$ be a positive solution of~$Q'[f]=0$ in~$\Omega$ and let~$w\in W^{1,p}(\Omega)\cap C_{c}(\Omega)$ be nonnegative. Write~$X(v,w)\triangleq\int_{\Omega}v^{p}|\nabla w|^{p}_{\mathcal{A}}\dx$ and~$X(w,v)\triangleq\int_{\Omega}w^{p}|\nabla v|^{p}_{\mathcal{A}}\dx$.
Then there exists a positive constant~$C$ (independent of~$v,w$) such that
\begin{align*}
Q[vw]\leq \begin{cases}		CX(v,w)&~\mbox{if}~p\leq 2,\\
		C\left(X(v,w)+X(v,w)^{2/p}X(w,v)^{1-2/p}\right)&~\mbox{if}~p>2.
		\end{cases}
\end{align*}
\end{corollary}
\subsection{Criticality theory}
\emph{Criticality theory} is the investigation of positivity properties of the operator $Q'$, which can be found in many mathematical works (see, e.g., \cite{HPR,PPAPDE,Regev} and references therein). For its applications, see, e.g., \cite{Giri1,Kovarik}.

In this subsection, we collect some results of criticality theory and review positive solutions of minimal growth.

The following Agmon--Allegretto--Piepenbrink-type (AAP-type) theorem, fundamental to criticality theory, connects the functional~$Q$ to the equation~$Q'[u]=0$.
\btheorem[{The AAP-type theorem \cite[Theorem~5.3]{HPR}}]\label{AAP}
The functional $Q$ is nonnegative in $\Gw$ if and only if the  equation $Q'[u]=0$ in $\Gw$ admits a positive (super-)solution in~$\Omega$.
\etheorem
  \btheorem[{\cite[Theorem~6.9]{HPR}}]\label{nullscon}
  Suppose that $Q$ is nonnegative in~$\Omega$. Then every null-sequence of~$Q$ in~$\Omega$ converges, both in~$L^{p}_{\loc}(\Omega)$ and a.e. in~$\Omega$, to a unique (up to a positive multiplicative constant) positive solution of the equation $Q'[u]=0$ in~$\Omega$. Up to a positive multiplicative constant, every ground state is the unique positive solution  and is also a unique positive supersolution (in~$W^{1,\infty}_{\loc}(\Omega)$ additionally for~$p<2$) of the same equation in~$\Omega$.
  \etheorem
By \cite[Theorem~6.12 (1)]{HPR}, the functional~$Q$ is critical in~$\Omega$ if and only if~$Q$ has a null-sequence in~$\core$. Checking its proof, we can see that the characterization (i) in the following theorem also holds.
  \begin{theorem}[{\cite[Theorem~6.12 (1) and (3)]{HPR}}]\label{nullc}
  Suppose that $Q$ is nonnegative in~$\Omega$.
  \begin{itemize}
      \item [\emph{(i)}]  The functional~$Q$ is critical in~$\Omega$ if and only if~$Q$ has a null-sequence.
      \item [\emph{(ii)}] The functional $Q$ is subcritical in $\Omega$ if and only if $Q$ admits a positive continuous Hardy-weight in $\Omega$.
  \end{itemize}
  \end{theorem}
 The following theorem consists of three perturbation results.
  \btheorem[{\cite[Proposition 6.14 (a) and Corollaries 6.15 and 6.17 ]{HPR}}]\label{perturbation}
  \leavevmode
\begin{itemize}
    \item [\emph{(i)}] Let $\Omega_{1}\subseteq\Omega_{2}$ be two subdomains of $\Omega$ with $\Omega_{2}\setminus\overline{\Omega_{1}}\neq\emptyset$. If $Q$ is nonnegative in $\Omega_{2}$, then $Q$ is subcritical in $\Omega_{1}$.
    \item [\emph{(ii)}] If $Q$ is nonnegative in $\Omega$, then the generalized principal eigenvalue $\lambda_{1}(Q;\omega)$ (see Definition~\ref{gpeigen}) is positive for all domains $\omega\Subset\Omega$.
    \item[\emph{(iii)}] For~$i=0,1$, let~$V_{i}\in\widehat{M}^{q}_{\loc}(p;\Omega)$. Assume that~$Q_{V_{i}}$ are nonnegative in~$\Omega$ for~$i=0,1$. Let~$V_{t}\triangleq(1-t)V_{0}+tV_{1} $ for~$t\in [0,1]$.
     Then $Q_{V_{t}}$ is nonnegative in~$\Omega$ for all~$t\in [0,1]$. Moreover, if $\vol\left(\{V_{0}\neq V_{1}\}\right)>0$, then~$Q_{V_{t}}$ is subcritical in~$\Omega$ for every~$t\in (0,1)$.
\end{itemize}
  \etheorem
  \bdefinition
\emph{A \emph{smooth (Lipschitz) exhaustion} of~$\Omega$ is a sequence of~$C^{\infty}$ (Lipschitz) domains~$\{\omega_{k}\}_{k\in\mathbb{N}}$ such that for all~$k\in\mathbb{N}$,~$\omega_{k}\Subset\omega_{k+1}\Subset\Omega$ and~$\cup_{k\in\mathbb{N}}\omega_{k}=\Omega$.}
\edefinition
\begin{definition}
\emph{A compact subset~$K$ of~$\Omega$ is called \emph{admissible} if~$K$ is the closure of a Lipschitz domain with connected boundary.}
\end{definition}
\begin{lemma}[{\cite[arXiv, Lemma~7.3]{HPR}}]
For every Lipschitz domain~$\omega_{0}\Subset\Omega$,~$K\triangleq\overline{\omega_{0}}$ is admissible if and only if for some (all) domain(s)~$\omega$ with~$K\subseteq\omega \Subset\Omega$,~$\omega\setminus K$ is a domain. In particular, for every admissible compact subset~$K$ of~$\Omega$ and  every Lipschitz exhaustion~$\{\omega_{i}\}_{i\in\mathbb{N}}$ of~$\Omega$, it holds that $\omega_{i}\setminus K$ is a domain for all sufficiently large~$i\in\mathbb{N}$.
\end{lemma}
\begin{remark}
\emph{By the above Lemma, it is easy to see that~$\Omega\setminus K$ is a domain for all admissible compact subsets~$K$ of~$\Omega$. See also \cite[Remark~4.22]{Hou2}.}
\end{remark}
We make some simple amendments to \cite[Definition 7.3 and arXiv, Definition 7.4]{HPR} and write down the following definition.          \begin{Def}\label{dfnmg}
 \emph{Let $K_{0}$ be a compact subset of~$\Omega$ such that~$\Omega\setminus K_{0}$ is a domain. A positive solution~$u$ of~$Q'[w]=0$ in~$\Omega\setminus K_{0}$ has \emph{minimal growth in a neighborhood of infinity} in $\Omega$ if for all admissible compact subsets~$K$ of~$\Omega$ with~$K_{0}\subseteq \mathring{K}$, and all positive solutions~$v\in C\left(\Omega\setminus \mathring{K}\right)$ of~$Q'[w]=g$  in~$\Omega\setminus K$ such that $u\leq v$ on~$\partial K$, it holds that~$ u\leq v$ in~$\Omega\setminus K$, where~$g\in M^{q}_{\loc}(p;\Omega)$ for $p\neq n$, $g\in L^{\bar{\rho}}_{\loc}(\Omega)$ ($\bar{\rho}>1$) for $p=n$, and $g$ is nonnegative a.e. in~$\Omega\setminus K$. The set of all such positive solutions is denoted by~$\mathcal{M}_{\Omega;K_{0}}$ or~$\mathcal{M}_{\Omega;K_{0}}^{V}$. If~$K_{0}=\emptyset$, then $u$ is called a \emph{global minimal positive solution} of~$Q'[w]=0$ in~$\Omega$.}
  \end{Def}
\section{Optimal Hardy-weights}\label{zpotential}
In this section, we introduce optimal Hardy-weights and prove some properties of them.

First, we define optimal Hardy-weights for the functional~$Q$. See also \cite[Definition~2.29]{Versano}. By Theorem \ref{nullscon} and Theorem \ref{nullc} (i), if the functional~$Q$ is critical, then up to a positive multiplicative constant, the equation~$Q'[u]=0$ has a unique positive solution, the functional~$Q$ has a unique ground state, and the ground state is the positive solution.
\bdefinition\label{ohw}
\emph{A nonnegative function~$W\in \widehat{M}^{q}_{\loc}(p;\Omega)$ 
 is called an \emph{optimal Hardy-weight} of~$Q$ in~$\Omega$ if
\begin{itemize}
\item (\textbf{Criticality}) the functional~$Q_{V-W}$ is critical in~$\Omega$; 
\item (\textbf{Null-criticality with respect to~$W$}) the ground state~$\psi$ of~$Q_{V-W}$ fulfills the condition~$\int_{\Omega}W\psi^{p}\dx=\infty$.
\end{itemize}
}
\edefinition
\bremarks\label{ohwdr}
\emph{\begin{itemize}
\item[\textup{(i)}] By the null-criticality condition, the zero function is never an optimal Hardy-weight. Furthermore, by the criticality condition, if~$Q$ has an optimal Hardy-weight in~$\Omega$, then~$Q$ is subcritical in~$\Omega$. 
\item[\textup{(ii)}] Given a subcritical functional~$Q$, for every nonnegative~$\phi\in\core\setminus\{0\}$, by \cite[Proposition~6.18]{HPR}, there is a positive constant~$\tau_{+}$ such that~$Q_{V-\tau_{+}\phi}$ is critical. However, for any such Hardy-weight~$\tau_{+}\phi$,~$Q_{V-\tau_{+}\phi}$ is not null-critical with respect to~$\tau_{+}\phi$. See also \cite[Remark~2.30]{Versano} for a similar discussion.
\end{itemize}}
\eremarks
By Theorem \ref{nullc} (ii), a subcritical functional~$Q$ has a positive continuous Hardy-weight in~$\Omega$.
\bdefinition\label{PAVHSS}
\emph{
Suppose that~$Q$ is subcritical in~$\Omega$ and let~$\mathcal{V}\in C(\Omega)\cap\mathcal{H}(\Omega)$ be positive in~$\Omega$. The \emph{$Q$-Sobolev space~$\widetilde{W}^{1,p}(\Omega)$} is defined as
\begin{align*}
\left\{\!u\!\in \!L^{p}(\Omega,(|V|+\mathcal{V})\dx)\bigg|\mbox{ the weak gradient }\nabla u\mbox{ exists in}~\Omega~\mbox{and}~\Vert u\Vert_{\widetilde{W}^{1,p}(\Omega)}\!<\!\infty\!\right\}\!,
\end{align*}
where
$$\Vert u\Vert_{\widetilde{W}^{1,p}(\Omega)}\triangleq Q_{|V|+\mathcal{V}}[u]^{1/p}\triangleq\left(\Vert|\nabla u|_{\mathcal{A}}\Vert^{p}_{L^{p}(\Omega)}+\Vert u\Vert^{p}_{L^{p}\left(\Omega,\left(|V|+\mathcal{V}\right)\dx\right)}\right)^{1/p}.$$
Let 
$\widetilde{W}^{1,p}_{0}(\Omega)=\widetilde{W}^{1,p}_{0,\mathcal{V}}(\Omega)$ be the closure of~$W^{1,p}(\Omega)\cap C_{c}(\Omega)$ in~$\widetilde{W}^{1,p}(\Omega)$.}
\edefinition
Checking \cite[p.~2]{Dintegral} and the proof of \cite[Lemma 4.5]{Hou2}, we can extract the following result.
\blemma[{\cite[Lemma 4.5]{Hou2}}]
Suppose that~$Q$ is subcritical in~$\Omega$. Then the~$Q$-Sobolev space~$\widetilde{W}^{1,p}(\Omega)$ is a Banach space and~$\widetilde{W}^{1,p}(\Omega)\subseteq W^{1,p}_{\loc}(\Omega)$.
\elemma
Recall that for every~$g\in\mathcal{H}(\Omega)$,~$S_{g}$ is the Hardy constant (see Definition \ref{hardyconstant}).
\blemma[{\cite[Lemma 4.9]{Hou2}}]\label{varhc}
Suppose that~$Q$ is subcritical in~$\Omega$. For every~$g\in \mathcal{H}(\Omega)$, 
$Q_{V-S_{g}|g|}[\phi]\geq 0$ for all~$\phi\in \widetilde{W}^{1,p}_{0}(\Omega)$ and$$S_{g}=\inf\bigg\{Q[\phi]~\bigg|~\phi\in \widetilde{W}^{1,p}_{0}(\Omega) \mbox{~and} \int_{\Omega}|g||\phi|^{p}\dx=1\bigg\}.$$
\elemma
The following corollary reveals that~$\widetilde{W}^{1,p}_{0,\mathcal{V}}(\Omega)$ is actually independent of~$\mathcal{V}$.
\bcorollary[{\cite[Corollary 4.10]{Hou2}}]\label{dndv}
Suppose that~$Q$ is subcritical in~$\Omega$ and let~$\mathcal{V},\mathcal{V}'\in C(\Omega)\cap\mathcal{H}(\Omega)$ be positive in~$\Omega$. Then~$\widetilde{W}^{1,p}_{0,\mathcal{V}}(\Omega)=\widetilde{W}^{1,p}_{0,\mathcal{V}'}(\Omega)$.
\ecorollary
\bremark[{\cite[Remark 4.11]{Hou2}}]\label{rem1}
\emph{The space~$\widetilde{W}^{1,p}_{0}(\Omega)$ is the completion of~$W^{1,p}(\Omega)\cap C_{c}(\Omega)$ with respect to the norm~$Q_{V^{+}}^{1/p}$. Therefore, when~$H(x,\xi)=\sqrt{A(x)\xi\cdot\xi}$ for almost all~$x\in\Omega$ and all~$\xi\in\R^{n}$, with the measurable matrix function~$A$ on~$\Omega$ satisfying some assumptions (see the conditions (1) and (2) in \cite[Theorem 3.22]{Hou2}), the space~$\widetilde{W}^{1,p}_{0}(\Omega)$ is equal to the \emph{generalized Beppo Levi space~$\mathcal{D}^{1,p}_{A,V^{+}}(\Omega)$} in \cite{Das} up to an isometry.}
\eremark
\bdefinition
\emph{Suppose that~$Q$ is subcritical in~$\Omega$. For every~$g\in\mathcal{H}(\Omega)\setminus\{0\}$, the Hardy constant~$S_{g}$ is \emph{attained at some~$\phi\in \widetilde{W}^{1,p}_{0}(\Omega)\setminus\{0\}$} if~$\int_{\Omega}|g||\phi|^{p}\dx>0$ and $$S_{g}=\frac{Q[\phi]}{\int_{\Omega}|g||\phi|^{p}\dx}.$$}
\edefinition
\blemma[{\cite[Lemma 4.13]{Hou2}}]\label{lemuni}
Suppose that~$Q$ is subcritical in~$\Omega$. Let~$g\in \mathcal{H}(\Omega)\cap M^{q}_{\loc}(p;\Omega)\setminus\{0\}$ and suppose that~$S_{g}$ is attained at some~$\phi\in \widetilde{W}^{1,p}_{0}(\Omega)\setminus\{0\}$. 
Then~$\phi$ is a solution of~$Q'_{V-S_{g}|g|}[u]=0$ in~$\Omega$ and either~$\phi$ or~$-\phi$ is positive in~$\Omega$. For every two such functions~$\phi$ and~$\phi'$, there exists~$C\in \R\setminus\{0\}$ such that~$\phi=C\phi'$.
\elemma
\blemma
Suppose that~$Q$ has an optimal Hardy-weight~$W$ in~$\Omega$. Then~$S_{W}=1$ and~$S_{W}$ is not attained at any function in~$\widetilde{W}^{1,p}_{0}(\Omega)\setminus\{0\}$.
\elemma
\bproof
Since~$Q_{V-W}$ is nonnegative in~$\Omega$, we get~$S_{W}\geq 1$. Clearly,~$S_{W}<\infty$. Because~$Q_{V-W}$ is critical in~$\Omega$,~$S_{W}>1$ is forbidden. Then~$S_{W}=1$.

Conversely, suppose that~$S_{W}$ is attained at some function $\phi\in \widetilde{W}^{1,p}_{0}(\Omega)\setminus\{0\}$. By Lemma \ref{lemuni},~$\phi$ is a solution of~$Q'_{V-W}[u]=0$ in~$\Omega$ and either~$\phi$ or~$-\phi$ is positive in~$\Omega$. We may assume that~$\phi$ is positive in~$\Omega$. Otherwise, consider~$-\phi$ because~$S_{W}$ is also attained at~$-\phi$ and~$-\phi$ is also a solution of~$Q'_{V-W}[u]=0$ in~$\Omega$. Then up to a positive multiplicative constant,~$\phi$ is the ground state of~$Q_{V-W}$. By the null-criticality with respect to~$W$,~$\int_{\Omega}W\phi^{p}\dx=\infty$. According to Lemma \ref{varhc}, we obtain~$Q[\phi]\geq\int_{\Omega}W\phi^{p}\dx$. Since~$\phi\in \widetilde{W}^{1,p}_{0}(\Omega)$,~$Q[\phi]<\infty$. Hence,~$\int_{\Omega}W\phi^{p}\dx<\infty$, a contradiction.
\eproof
\blemma
Suppose that~$Q$ has an optimal Hardy-weight~$W$ in~$\Omega$. If~$V_{1}\in \widehat{M}^{q}_{\loc}(p;\Omega)$ fulfills~$V_{1}\geq -\varepsilon W$ for some~$\varepsilon<1$, then~$W+V_{1}$ is an optimal Hardy-weight of~$Q_{V+V_{1}}$ in~$\Omega$.
\elemma
\bproof
The proof is similar to that of \cite[Lemma~2.35]{Versano}. Note that~$W+V_{1}\in\widehat{M}^{q}_{\loc}(p;\Omega)$ and that~$Q_{V+V_{1}-(W+V_{1})}=Q_{V-W}$ is critical in~$\Omega$.  In addition,~$W+V_{1}\geq (1-\varepsilon)W\geq0$. Let~$\psi$ be the ground state of~$Q_{V-W}$ in~$\Omega$. Then~$\psi$ is also the ground state of~$Q_{V+V_{1}-(W+V_{1})}$ in~$\Omega$. We also have 
$$\int_{\Omega}(W+V_{1})\psi^{p}\dx\geq (1-\varepsilon)\int_{\Omega}W\psi^{p}\dx=\infty.$$ Thus~$W+V_{1}$ is an optimal Hardy-weight of~$Q_{V+V_{1}}$ in~$\Omega$. 
\eproof
The following lemma is an analogue of \cite[Theorem~3.1]{Kovarik}. The proof is similar and hence omitted.
\blemma\label{lem216kp}
Suppose that~$Q$ is critical in~$\Omega$ and that~$Q_{V-W}$ is nonnegative in~$\Omega\setminus\overline{\omega}$ for some nonnegative~$W\in L^{1}_{\loc}(\Omega)$ and some domain~$\omega\Subset\Omega$. Let~$\phi$ be the ground state of~$Q$ in~$\Omega$. Then for all domains~$\omega'\Subset\Omega$ such that~$\omega\Subset\omega'$, it holds that~$\int_{\Omega\setminus\overline{\omega'}}W|\phi|^{p}\dx<\infty$.
\elemma
\blemma\label{optimalinf}
If~$W$ is an optimal Hardy-weight of~$Q$ in~$\Omega$, then for all compact sets~$K\subseteq\Omega$, 
$$\sup\{\lambda\in\R~|~Q_{V-\lambda W}\geq 0~\mbox{in}~ \Omega\setminus K\}=\max\{\lambda\in\R~|~Q_{V-\lambda W}\geq 0~\mbox{in}~ \Omega\setminus K\}=1.$$ 
\elemma
\bproof
The proof is similar to that of \cite[Corollary~3.4]{Kovarik}. To obtain a contradiction, suppose that there exists a compact set~$K\subseteq\Omega$ and a constant~$\lambda>1$ such that~$Q_{V-\lambda W}$ is nonnegative in~$\Omega\setminus K$. Therefore,~$Q_{V-W-(\lambda-1)W}$ is nonnegative in~$\Omega\setminus \overline{\omega}$ for some domain~$\omega\Subset\Omega$ such that~$K\subseteq\omega$. Let~$\phi$ be the ground state of~$Q_{V-W}$ in~$\Omega$. By Lemma~\ref{lem216kp}, for some domain~$\omega'\Subset\Omega$, it holds that~$\int_{\Omega\setminus\overline{\omega'}}W|\phi|^{p}\dx<\infty$. 
It follows that~$\int_{\Omega}W|\phi|^{p}\dx<\infty$, contrary to the null-criticality with respect to~$W$.
\eproof
\section{Optimal Hardy-weights: Zero potential}\label{zeropcons}
In this section, we first prove a corea-type formula and then construct optimal Hardy-weights for the functional~$Q_{0}[\cdot;\Omega^{*}]$. We also give some examples of Finsler~$p$-harmonic functions realizing relevant boundary conditions.
\bdefinition\label{aharmonic}
\emph{A $Q_{0}$-harmonic function in~$\Omega$ is also called~\emph{Finsler $p$-harmonic} in~$\Omega$.
}
\edefinition
We recall the definition of proper maps.
\bdefinition
\emph{Let~$X$ and~$Y$ be two topological spaces. A map~$f:X\rightarrow Y$ is \emph{proper} if~$f^{-1}(K)$ is compact for every compact subset~$K$ of~$Y$.}
\edefinition
The following coarea-type lemma will be used exclusively in Theorem~\ref{zerop}. In this lemma, it is sufficient that Assumptions~\ref{ass9} and~\ref{C1alpha} hold in~$\Omega^{*}$.
\blemma\label{coarea}
Let~$G\in C^{1}(\Omega^{*})$ be a positive Finsler $p$-harmonic function in~$\Omega^{*}$ such that~$G(\Omega^{*})=(\inf_{\Omega^{*}}G,\sup_{\Omega^{*}}G)$ and let~$h\in C^{2}((\inf_{\Omega^{*}}G,\sup_{\Omega^{*}}G))$ be a positive function. Suppose that~$G:\Omega^{*}\rightarrow (\inf_{\Omega^{*}}G,\sup_{\Omega^{*}}G)$ is proper. Let~$v\triangleq h(G)$. Then with respect to almost all~$\inf_{\Omega^{*}}G<t<\sup_{\Omega^{*}}G$,~$\int_{G^{-1}(t)}\frac{|\nabla G|_{\mathcal{A}}^{p}}{|\nabla G|}\dHnn$ is a positive constant, denoted by~$C$, and
$$\int_{\Omega^{*}}f(v)|\nabla v|_{\mathcal{A}}^{p}\dx=C\int_{\inf_{\Omega^{*}}G}^{\sup_{\Omega^{*}}G}f(h(t))|h'(t)|^{p}\dt,$$
where $f$ is a function on~$(0,\infty)$ such that~$f(v)\in L^{1}_{\loc}(\Omega^{*})$ is compactly supported in~$\Omega^{*}$ and $f(h(t))$ is measurable on~$(\inf_{\Omega^{*}}G,\sup_{\Omega^{*}}G)$.
\elemma
\bproof
We first note that for almost all~$\inf_{\Omega^{*}} G<t_{1}<t_{2}<\sup_{\Omega^{*}}G$,~$G^{-1}((t_{1},t_{2}))$ is bounded and~$\partial G^{-1}((t_{1},t_{2}))=G^{-1}(t_{1})\cup G^{-1}(t_{2})$. Since~$G:\Omega^{*}\rightarrow (\inf_{\Omega^{*}}G,\sup_{\Omega^{*}}G)$ is proper,~$G^{-1}([t_{1},t_{2}])$ is compact for all~$\inf_{\Omega^{*}} G<t_{1}<t_{2}<\sup_{\Omega^{*}}G$. Therefore,~$G^{-1}((t_{1},t_{2}))$ is bounded and~$\partial G^{-1}((t_{1},t_{2}))\subseteq\overline{G^{-1}((t_{1},t_{2}))}\subseteq G^{-1}([t_{1},t_{2}])\subseteq\Omega^{*}$ for all~$\inf_{\Omega^{*}} G<t_{1}<t_{2}<\sup_{\Omega^{*}}G$. Since~$G$ is continuous in~$\Omega^{*}$,~$\partial G^{-1}((t_{1},t_{2}))\subseteq G^{-1}(t_{1})\cup G^{-1}(t_{2})$ for all~$\inf_{\Omega^{*}} G<t_{1}<t_{2}<\sup_{\Omega^{*}}G$. By the generalized Sard theorem \cite[Theorem~1.2]{Sard}, for almost all~$\inf_{\Omega^{*}}G<t<\sup_{\Omega^{*}}G$,$$\mathcal{H}^{n-1}\left(\{x\in\Omega^{*}~|~G(x)=t~\mbox{and}~\nabla G(x)=0\}\right)=0.$$ By virtue of the strong maximum principle \cite[6.5]{HKM}, we deduce that~$\partial G^{-1}((t_{1},t_{2}))\supseteq G^{-1}(t_{1})\cup G^{-1}(t_{2})$ for almost all~$\inf_{\Omega^{*}} G<t_{1}<t_{2}<\sup_{\Omega^{*}}G$. Then~$\partial G^{-1}((t_{1},t_{2}))=G^{-1}(t_{1})\cup G^{-1}(t_{2})$ for almost all~$\inf_{\Omega^{*}} G<t_{1}<t_{2}<\sup_{\Omega^{*}}G$.

The remaining proof is similar to that of \cite[Lemma~2.33]{Versano}. Suppose that~$\supp f(v)\subseteq\omega_{0}$ for some domain~$\omega_{0}\Subset\Omega^{*}$.
By \cite[Theorem~2.32]{ChenCPAM}, we may deduce that
\begin{align}\label{fcalc}
&\int_{\Omega^{*}}f(v)|\nabla v|_{\mathcal{A}}^{p}\dx= \int_{\Omega^{*}}f(h(G))|h'(G)|^{p}|\nabla G|_{\mathcal{A}}^{p}\dx=\int_{\omega_{0}}f(h(G))|h'(G)|^{p}|\nabla G|_{\mathcal{A}}^{p}\dx\notag\\
&=\int_{\omega_{0}}f(h(G))|h'(G)|^{p}\frac{|\nabla G|_{\mathcal{A}}^{p}}{|\nabla G|}|\nabla G|\dx\notag\\
&=\int_{\inf_{\Omega^{*}}G}^{\sup_{\Omega^{*}}G}\int_{\{x\in\omega_{0}|\,G(x)=t\}}f(h(G(x)))|h'(G(x))|^{p}\frac{|\nabla G|_{\mathcal{A}}^{p}}{|\nabla G|}\dHnn\dt\notag\\
&=\int_{\inf_{\Omega^{*}}G}^{\sup_{\Omega^{*}}G}\int_{\{x\in\Omega^{*}|\,G(x)=t\}}f(h(G(x)))|h'(G(x))|^{p}\frac{|\nabla G|_{\mathcal{A}}^{p}}{|\nabla G|}\dHnn\dt\notag\\
&=\int_{\inf_{\Omega^{*}}G}^{\sup_{\Omega^{*}}G}f(h(t))|h'(t)|^{p}\int_{G^{-1}(t)}\frac{|\nabla G|_{\mathcal{A}}^{p}}{|\nabla G|}\dHnn\dt.
\end{align}
Note that for almost all~$\inf_{\Omega^{*}}G<t_{1}<t_{2}<\sup_{\Omega^{*}}G$,~$G^{-1}((t_{1},t_{2}))$ is compactly contained in~$\Omega^{*}$ and is an open set with almost~$C^{1}$-boundary (see \cite[Section~9.3]{Maggi}). By virtue of the Gauss--Green theorem on open sets with almost~$C^{1}$-boundary \cite[Theorem~9.6]{Maggi}, we may conclude that~$G^{-1}((t_{1},t_{2}))$ has finite perimeter in~$\Omega^{*}$. Recall (see Lemma \ref{pflemma}) that for almost all~$x\in\Omega^{*}$ and all~$\xi\in\R^{n}$,$$\mathcal{A}(x,\xi)\cdot\xi=pF(x,\xi)=H(x,\xi)^{p}=|\xi|_{\mathcal{A}}^{p}.$$ Since~$-\dive\mathcal{A}(x,\nabla G)=0$ in~$\Omega^{*}$, by \cite[Theorems~5.2 and 7.2]{ChenCPAM}, we may calculate: for almost all~$\inf_{\Omega^{*}}G<t_{1}<t_{2}<\sup_{\Omega^{*}}G$,
\begin{align*}
0&=\int_{G^{-1}(t_{2})}\mathcal{A}(x,\nabla G)\cdot\frac{\nabla G}{|\nabla G|}\dHnn-\int_{G^{-1}(t_{1})}\mathcal{A}(x,\nabla G)\cdot\frac{\nabla G}{|\nabla G|}\dHnn\\
&=\int_{G^{-1}(t_{2})}\frac{|\nabla G|_{\mathcal{A}}^{p}}{|\nabla G|}\dHnn-\int_{G^{-1}(t_{1})}\frac{|\nabla G|_{\mathcal{A}}^{p}}{|\nabla G|}\dHnn.
\end{align*}
Therefore, with respect to almost all~$\inf_{\Omega^{*}}G<t<\sup_{\Omega^{*}}G$,~$\int_{G^{-1}(t)}\frac{|\nabla G|_{\mathcal{A}}^{p}}{|\nabla G|}\dHnn$ is a positive (see \cite[Theorem 3.10]{evansf}) constant, denoted by~$C$. Then by \eqref{fcalc},
\begin{align*}
\int_{\Omega^{*}}f(v)|\nabla v|_{\mathcal{A}}^{p}\dx=C\int_{\inf_{\Omega^{*}}G}^{\sup_{\Omega^{*}}G}f(h(t))|h'(t)|^{p}\dt.\qquad\qedhere
\end{align*}
\eproof
\bremark\label{propernece}
\emph{In this lemma, under all the other conditions, if for almost all~$\inf_{\Omega^{*}} G<t_{1}<t_{2}<\sup_{\Omega^{*}}G$,~$G^{-1}((t_{1},t_{2}))$ is bounded and~$\partial G^{-1}((t_{1},t_{2}))=G^{-1}(t_{1})\cup G^{-1}(t_{2})$, then~$G:\Omega^{*}\rightarrow (\inf_{\Omega^{*}}G,\sup_{\Omega^{*}}G)$ is proper. The proof is elementary and hence omitted.
}
\eremark
The proof of the following lemma is similar to those of \cite[arXiv, Lemma~A.1]{Das} and \cite[Lemma~2.10]{Devyver1}. See also \cite[Proposition~3.4]{Versano}.
\blemma\label{chlem}
Suppose that~$f\in C^{2}((0,\infty))$ and that~$f$ and~$f'$ are both positive on~$(0,\infty)$. Then for all positive~$u\in C^{1}(\Omega)$ and all~$\phi\in\core$,
\begin{align*}
&\int_{\Omega}\mathcal{A}(x,\nabla (f(u)))\cdot\nabla\phi\dx\\
&=-\int_{\Omega}(p-1)f'(u)^{p-2}f''(u)|\nabla u|_{\mathcal{A}}^{p}\phi\dx+\int_{\Omega}\mathcal{A}(x,\nabla u)\cdot\nabla(f'(u)^{p-1}\phi)\dx,
\end{align*}
{and in the weak sense,
\begin{align*}
-\dive\mathcal{A}(x,\nabla (f(u)))=-(p-1)f'(u)^{p-2}f''(u)|\nabla u|_{\mathcal{A}}^{p}-f'(u)^{p-1}\dive\mathcal{A}(x,\nabla u).
\end{align*}}
\elemma
\bproof
By the product and chain rules, we may write:
\begin{align*}
&\int_{\Omega}\mathcal{A}(x,\nabla (f(u)))\cdot\nabla\phi\dx\\
&=-\int_{\Omega}\mathcal{A}(x,\nabla u)\cdot\nabla(f'(u)^{p-1})\phi\dx+\int_{\Omega}\mathcal{A}(x,\nabla u)\cdot\nabla(f'(u)^{p-1}\phi)\dx\\
&=-\int_{\Omega}(p-1)f'(u)^{p-2}f''(u)|\nabla u|_{\mathcal{A}}^{p}\phi\dx+\int_{\Omega}\mathcal{A}(x,\nabla u)\cdot\nabla(f'(u)^{p-1}\phi)\dx.\qedhere
\end{align*}
\eproof
Recall that~$f_{0}(t)=t^{(p-1)/p}$ for all~$t>0$. Then
$$f_{0}\in C^{\infty}(0,\infty),\quad f'_{0}(t)=\frac{p-1}{p}t^{-1/p},\quad\mbox{and}\quad f''_{0}(t)=-\frac{p-1}{p^{2}}t^{-1/p-1}.$$
The following corollary is a generalization of \cite[Corollary~3.6]{Versano} and the proof is similar.
\bcorollary\label{co2}
Suppose that~$G_{\phi}$ is a positive solution of~$Q'_{c_{p}V}[u]=\phi$ in~$\Omega$, where~$0\leq\phi\in\core$ and~$c_{p}=(p/(p-1))^{p-1}$. Then
 for all~$\varphi\in\core$, 
\begin{align*}
&\int_{\Omega}\mathcal{A}(x,\nabla (f_{0}(G_{\phi})))\cdot\nabla\varphi\dx+\int_{\Omega}Vf_{0}(G_{\phi})^{p-1}\varphi\dx\\
&=-\int_{\Omega}(p-1)f'_{0}(G_{\phi})^{p-2}f_{0}''(G_{\phi})|\nabla G_{\phi}|_{\mathcal{A}}^{p}\varphi\dx+\int_{\Omega}f'_{0}(G_{\phi})^{p-1}\phi\varphi\dx\\
&=\int_{\Omega}\left(\frac{p-1}{p}\right)^{p}\left|\frac{\nabla G_{\phi}}{G_{\phi}}\right|^{p}_{\mathcal{A}}f_{0}(G_{\phi})^{p-1}\varphi\dx+\left(\frac{p-1}{p}\right)^{p-1}\int_{\Omega}G_{\phi}^{(1-p)/p}\phi\varphi\dx
\end{align*}
\ecorollary
\bproof
By Lemma~\ref{chlem}, for all~$\varphi\in\core$, 
\begin{align*}
&\int_{\Omega}\mathcal{A}(x,\nabla (f_{0}(G_{\phi})))\cdot\nabla\varphi\dx\\
&=-\int_{\Omega}(p-1)f'_{0}(G_{\phi})^{p-2}f''_{0}(G_{\phi})|\nabla G_{\phi}|_{\mathcal{A}}^{p}\varphi\dx+\int_{\Omega}\mathcal{A}(x,\nabla G_{\phi})\cdot\nabla(f_{0}'(G_{\phi})^{p-1}\varphi)\dx
\end{align*}
Since~$f_{0}'(G_{\phi})^{p-1}\varphi\in C^{1}_{c}(\Omega)\subseteq W^{1,p}_{c}(\Omega)$ and~$Q'_{c_{p}V}[G_{\phi}]=\phi$, we deduce that
$$\int_{\Omega}\mathcal{A}(x,\nabla G_{\phi})\cdot\nabla(f_{0}'(G_{\phi})^{p-1}\varphi)\dx+\int_{\Omega}c_{p}VG_{\phi}^{p-1}f_{0}'(G_{\phi})^{p-1}\varphi\dx=\int_{\Omega}\phi f_{0}'(G_{\phi})^{p-1}\varphi\dx.$$
Note that for all~$t>0$,
$$c_{p}t^{p-1}f'_{0}(t)^{p-1}=\left(\frac{p}{p-1}\right)^{p-1}t^{p-1}\left(\frac{p-1}{p}t^{-1/p}\right)^{p-1}=t^{p-2+1/p}=t^{(p-1)^{2}/p}=f_{0}(t)^{p-1}.$$
Therefore, $$\int_{\Omega}Vf_{0}(G_{\phi})^{p-1}\varphi\dx=\int_{\Omega}c_{p}VG_{\phi}^{p-1}f_{0}'(G_{\phi})^{p-1}\varphi\dx.$$
Furthermore, for all~$\varphi\in\core$, 
\begin{align*}
&\int_{\Omega}\mathcal{A}(x,\nabla (f_{0}(G_{\phi})))\cdot\nabla\varphi\dx+\int_{\Omega}Vf_{0}(G_{\phi})^{p-1}\varphi\dx\\
&=-\int_{\Omega}(p-1)f'_{0}(G_{\phi})^{p-2}f''_{0}(G_{\phi})|\nabla G_{\phi}|_{\mathcal{A}}^{p}\varphi\dx+\int_{\Omega}\mathcal{A}(x,\nabla G_{\phi})\cdot\nabla(f_{0}'(G_{\phi})^{p-1}\varphi)\dx\\
&\hspace{5mm}+\int_{\Omega}c_{p}VG_{\phi}^{p-1}f_{0}'(G_{\phi})^{p-1}\varphi\dx
\\
&=-\int_{\Omega}(p-1)f'_{0}(G_{\phi})^{p-2}f''_{0}(G_{\phi})|\nabla G_{\phi}|_{\mathcal{A}}^{p}\varphi\dx+\int_{\Omega}f_{0}'(G_{\phi})^{p-1}\phi\varphi\dx\\
&=\int_{\Omega}\left(\frac{p-1}{p}\right)^{p}\left|\frac{\nabla G_{\phi}}{G_{\phi}}\right|^{p}_{\mathcal{A}}f_{0}(G_{\phi})^{p-1}\varphi\dx+\left(\frac{p-1}{p}\right)^{p-1}\int_{\Omega}G_{\phi}^{(1-p)/p}\phi\varphi\dx.
\qedhere
\end{align*}
\eproof
In the following theorem, it is also sufficient that Assumptions~\ref{ass9} and~\ref{C1alpha} hold in~$\Omega^{*}$. Recall that~$f_{0}(t)=t^{(p-1)/p}$ for all~$t>0$.
\begin{theorem}\label{zerop}
Suppose that Assumption~\ref{assup} holds in~$\Omega^{*}$ and that there exists a positive Finsler $p$-harmonic function~$\mathcal{G}$ in~$\Omega^{*}$ satisfying, for some nonnegative real number~$\sigma$,
 \begin{align}\label{bcons}
\begin{cases}		\lim_{x\rightarrow\hat{x}}\mathcal{G}(x)=\infty\quad\mbox{and}\quad\lim_{x\rightarrow \overline{\infty}}\mathcal{G}(x)=0&\mbox{if}~1<p\leq n,\\
\lim_{x\rightarrow\hat{x}}\mathcal{G}(x)=\sigma\quad\mbox{and}\quad\lim_{x\rightarrow \overline{\infty}}\mathcal{G}(x)=\begin{cases}
    \infty&~ \mbox{if}~\sigma=0,\\    0&~ \mbox{if}~\sigma>0,
\end{cases}&\mbox{if}~p>n.
		\end{cases}
\end{align}
Let~$W_{0}\triangleq\left(\frac{p-1}{p}\right)^{p}\left|\frac{\nabla\mathcal{G}}{\mathcal{G}}\right|_{\mathcal{A}}^{p}$ and
\begin{align*}
\begin{cases} 
v\triangleq f_{0}(\mathcal{G})~\mbox{and}~W\triangleq W_{0}&\mbox{if}~1<p\leq n, \mbox{or}~p>n~\mbox{and}~\sigma=0,\\ v\triangleq f_{0}(\mathcal{G}(\sigma-\mathcal{G}))~\mbox{and}~W\triangleq\\\left(\frac{1}{\sigma-\mathcal{G}}\right)^{p}W_{0}|\sigma-2\mathcal{G}|^{p-2}(2(p-2)\mathcal{G}(\sigma-\mathcal{G})+\sigma^{2})&\mbox{otherwise}.
\end{cases}
\end{align*}
Then in~$\Omega^{*}$,~$W$ is an optimal Hardy-weight of~$Q_{0}$ and (up to a positive multiplicative constant)~$v$ is the unique ground state of~$Q_{-W}$. 
\end{theorem}
\bproof
We may assume that~$\hat{x}=0$. Otherwise, consider the translation~$\Omega-\hat{x}$.

The proof is similar to those of \cite[Propositions 5.1, 5.2, 5.5, and 5.6]{Devyver1}. 

We first demonstrate that~$\mathcal{G}<\sigma$ in~$\Omega^{*}$ if~$p>n$ and~$\sigma>0$. Let~$\{\omega_{i}\}_{i\in\mathbb{N}}$ be a smooth exhaustion of~$\Omega$. Since $\lim_{x\rightarrow 0}\mathcal{G}(x)=\sigma$, for each~$\varepsilon>0$, there exists~$\delta>0$ such that~$\mathcal{G}(x)<\sigma+\varepsilon$ on~$\overline{B_{\delta}(0)}\setminus\{0\}$. In addition, since $\lim_{x\rightarrow \overline{\infty}}\mathcal{G}(x)=0$, there exists~$N\in\mathbb{N}$ such that~$\mathcal{G}<\sigma$ on~$\partial\omega_{i}$ for all~$i>N$. We may also 
 assume that~$B_{\delta}(0)\Subset\omega_{i}$ for all~$i>N$. By the weak comparison principle \cite[Theorem~4.25]{HPR},~$\mathcal{G}\leq \sigma+\varepsilon$ on~$\overline{\omega_{i}}\setminus B_{\delta'}(0)$ for all~$i>N$ and all~$0<\delta'\leq\delta$. Then~$\mathcal{G}\leq \sigma+\varepsilon$ in~$\Omega^{*}$. Furthermore,~$\mathcal{G}\leq \sigma$ in~$\Omega^{*}$ because~$\varepsilon>0$ is arbitrary. By the strong maximum principle \cite[6.5]{HKM}, either~$\mathcal{G}<\sigma$ in~$\Omega^{*}$ or~$\mathcal{G}=\sigma$ in~$\Omega^{*}$. Because $\lim_{x\rightarrow \overline{\infty}}\mathcal{G}(x)=0$, we deduce that~$\mathcal{G}<\sigma$ in~$\Omega^{*}$. 

By virtue of \eqref{bcons}, in all the three cases, it is easy to show that\begin{align*}
\mathcal{G}(\Omega^{*})=\left(\inf_{\Omega^{*}}\mathcal{G},\sup_{\Omega^{*}}\mathcal{G}\right)=\begin{cases} 
(0,\infty)&\mbox{if}~1<p\leq n, \mbox{or}~p>n~\mbox{and}~\sigma=0,\\ (0,\sigma)&\mbox{otherwise},
\end{cases}
\end{align*} and that~$\mathcal{G}:\Omega^{*}\rightarrow(\inf_{\Omega^{*}}\mathcal{G},\sup_{\Omega^{*}}\mathcal{G})$ is proper.

 It is clear that~$W$ is nonnegative in~$\Omega^{*}$ in both cases. Since~$\mathcal{G}\in C^{1}(\Omega^{*})$,~$\mathcal{G}>0$ in~$\Omega^{*}$, and in the second case,~$\mathcal{G}<\sigma$ in~$\Omega^{*}$, it is easily seen that~$W\in L^{\infty}_{\loc}(\Omega^{*})\subseteq\widehat{M}^{q}_{\loc}(p;\Omega^{*})$ in both cases. 

Next we prove the criticality of~$Q_{-W}$ in~$\Omega^{*}$, which falls naturally into two parts.
\begin{itemize}
\item Suppose that~$1<p\leq n$, or~$p>n$ and~$\sigma=0$. By Corollary~\ref{co2},~$v$ is a positive solution of the equation~$Q'_{-W}[u]=0$ in~$\Omega^{*}$. By the AAP-type theorem (Theorem~\ref{AAP}),~$Q_{-W}$ is nonnegative in~$\Omega^{*}$. We next aim to find a null-sequence with respect to~$Q_{-W}$ in~$\Omega^{*}$ so that by Theorem~\ref{nullc}~(i),~$Q_{-W}$ is critical in~$\Omega^{*}$. For all~$k\in\mathbb{N}$ ($k\geq 2$), let
\begin{align*}
\varphi_{k}(t)\triangleq\begin{cases}
				0&0\leq t\leq \frac{1}{k^{2}},\\
			2+\frac{\log t}{\log k}&\frac{1}{k^{2}}\leq t\leq \frac{1}{k},\\
  1&\frac{1}{k}\leq t\leq k,\\
  2-\frac{\log t}{\log k}&k\leq t\leq k^{2},\\
  0& t\geq k^{2},
		\end{cases}
\end{align*}
$w_{k}\triangleq\varphi_{k}(v)$, and~$u_{k}\triangleq vw_{k}$. {Evidently,~$v^{-1}((1,2))$ is open,~$\emptyset\neq v^{-1}((1,2))\Subset\Omega^{*}$, and for all~$k\in\mathbb{N}$ ($k\geq 2$),~$\int_{v^{-1}((1,\,2))}u_{k}^{p}\dx\geq\vol(v^{-1}((1,2)))>0$. We claim that~$\{u_{k}/\Vert u_{k}\Vert_{L^{p}(v^{-1}((1,\,2)))}\}_{k\in\mathbb{N}\,(k\geq 2)}$ is a desired null-sequence.} Note that~$\varphi_{k}$ is Lipschitz continuous on~$[0,\infty)$ for all~$k\in\mathbb{N}$ ($k\geq 2$). Then~$w_{k}$ is locally Lipschitz continuous in~$\Omega^{*}$ and hence belongs to~$W^{1,\infty}_{\loc}(\Omega^{*})$ for all~$k\in\mathbb{N}$ ($k\geq 2$). It can be easily seen that~$w_{k}$ has compact support in~$\Omega^{*}$ for all~$k\in\mathbb{N}$ ($k\geq 2$) because of the relevant conditions in \eqref{bcons}. Thus~$w_{k}\in W^{1,p}(\Omega^{*})\cap C_{c}(\Omega^{*})$ for all~$k\in\mathbb{N}$ ($k\geq 2$). Furthermore,~$u_{k}\in W^{1,p}(\Omega^{*})\cap C_{c}(\Omega^{*})$ for all~$k\in\mathbb{N}$ ($k\geq 2$) because~$v\in C^{1}(\Omega^{*})$. Obviously,~$w_{k}$ and~$u_{k}$ are nonnegative in~$\Omega^{*}$ for all~$k\in\mathbb{N}$ ($k\geq 2$). By Corollary~\ref{co1}, there exists a positive constant~$C$, independent of~$k$, such that for all~$k\in\mathbb{N}$ ($k\geq 2$),
\begin{align}\label{se1}
Q_{-W}[u_{k};\Omega^{*}]\leq \begin{cases}		CX(v,w_{k})&~\mbox{if}~p\leq 2,\\
		C\left(X(v,w_{k})+X(v,w_{k})^{2/p}X(w_{k},v)^{1-2/p}\right)&~\mbox{if}~p>2,
		\end{cases}
\end{align}
where~$X(v,w_{k})=\int_{\Omega^{*}}v^{p}|\nabla w_{k}|^{p}_{\mathcal{A}}\dx$ and~$X(w_{k},v)=\int_{\Omega^{*}}w_{k}^{p}|\nabla v|^{p}_{\mathcal{A}}\dx$. By Lemma~\ref{coarea} with $h=f_{0}$ and $f(t)=t^{p}|\varphi_{k}'(t)|^{p}$ for all $t>0$, there exists a positive constant~$C'$ such that for all~$k\in\mathbb{N}$ ($k\geq 2$),
\begin{align*}
&X(v,w_{k})=\int_{\Omega^{*}}v^{p}|\nabla w_{k}|^{p}_{\mathcal{A}}\dx=\int_{\Omega^{*}}v^{p}|\varphi_{k}'(v)|^{p}|\nabla v|^{p}_{\mathcal{A}}\dx\\
&=C'\int_{0}^{\infty}f_{0}(t)^{p}|\varphi_{k}'(f_{0}(t))|^{p}f_{0}'(t)^{p}\dt=C'C(p)\int_{0}^{\infty}|\varphi_{k}'(f_{0}(t))|^{p}t^{p-2}\dt\\
&=C'C(p)\left(\frac{1}{\log k}\right)^{p}\left(\int_{(1/k^{2})^{p/(p-1)}}^{(1/k)^{p/(p-1)}}\frac{1}{t}\dt+\int_{k^{p/(p-1)}}^{k^{2p/(p-1)}}\frac{1}{t}\dt\right)\\
&=C'C(p)\left(\frac{1}{\log k}\right)^{p-1}
\end{align*}
and \begin{align*}
&X(w_{k},v)=\int_{\Omega^{*}}w_{k}^{p}|\nabla
v|^{p}_{\mathcal{A}}\dx=\int_{\Omega^{*}}\varphi_{k}(v)^{p}|\nabla v|^{p}_{\mathcal{A}}\dx=C'\int_{0}^{\infty}\varphi_{k}(f_{0}(t))^{p}f'_{0}(t)^{p}\dt\\&=C'C(p)\int_{0}^{\infty}\varphi_{k}(f_{0}(t))^{p}t^{-1}\dt\leq C'C(p)\int_{(1/k^{2})^{p/(p-1)}}^{k^{2p/(p-1)}}\frac{1}{t}\dt=C'C(p)\log k,
\end{align*}
where we stipulate that~$\varphi_{k}'(1/k^{2})=\varphi_{k}'(1/k)=\varphi_{k}'(k)=\varphi_{k}'(k^{2})=0$. By \eqref{se1}, we obtain~$\lim_{k\rightarrow\infty}Q_{-W}[u_{k};\Omega^{*}]=0$. {Then the claim follows.}

\item Suppose that~$p>n$ and~$\sigma>0$. By considering~$\mathcal{G}/\sigma$, without loss of generality, we can suppose that~$\sigma=1$. {Let~$g(s)\triangleq f_{0}(s(1-s))=(s(1-s))^{(p-1)/p}$ for all~$s\in (0,1)$}. Then~$g\in C^{2}(0,1)$ and $v=g(\mathcal{G})$. Since~$p>n\geq 2$,~$\varphi(s)\triangleq|s|^{p-2}s$ is continuously differentiable on~$\R$ and that~$\varphi'(s)=(p-1)|s|^{p-2}$. Note that in~$\Omega^{*}$,
 \begin{align*}
 &-(p-1)|g'(\mathcal{G})|^{p-2}g''(\mathcal{G})|\nabla \mathcal{G}|_{\mathcal{A}}^{p}\\
 =&(p-1)\frac{p-1}{p}|\nabla \mathcal{G}|_{\mathcal{A}}^{p}\left(\frac{p-1}{p}(\mathcal{G}(1-\mathcal{G}))^{-1/p}|1-2\mathcal{G}|\right)^{p-2}\times\\
 &\left(\frac{1}{p}(\mathcal{G}(1-\mathcal{G}))^{-1/p-1}(1-2\mathcal{G})^{2}+2(\mathcal{G}(1-\mathcal{G}))^{-1/p}\right)\\
 =&\left(\frac{p-1}{p}\right)^{p}|\nabla\mathcal{G}|_{\mathcal{A}}^{p}\left((\mathcal{G}(1-\mathcal{G}))^{-1/p}|1-2\mathcal{G}|\right)^{p-2}\times\\
 &\left((\mathcal{G}(1-\mathcal{G}))^{-1/p-1}(1-2\mathcal{G})^{2}+2p(\mathcal{G}(1-\mathcal{G}))^{-1/p}\right)\\
 =&\left(\frac{p-1}{p}\right)^{p}\left|\frac{\nabla\mathcal{G}}{\mathcal{G}}\right|_{\mathcal{A}}^{p}\left(\frac{1}{1-\mathcal{G}}\right)^{p}|1-2\mathcal{G}|^{p-2}(2(p-2)\mathcal{G}(1-\mathcal{G})+1)(\mathcal{G}(1-\mathcal{G}))^{(p-1)^{2}/p}\\
 =&Wv^{p-1}.
 \end{align*}Then as in Lemma~\ref{chlem}, for all~$\phi\in C^{\infty}_{c}(\Omega^{*})$,\begin{align*}
&\int_{\Omega^{*}}\mathcal{A}(x,\nabla v)\cdot\nabla\phi\dx\\
&=-\int_{\Omega^{*}}\mathcal{A}(x,\nabla \mathcal{G})\cdot\nabla(|g'(\mathcal{G})|^{p-2}g'(\mathcal{G}))\phi\dx\\
&~~~~+\int_{\Omega^{*}}\mathcal{A}(x,\nabla \mathcal{G})\cdot\nabla(|g'(\mathcal{G})|^{p-2}g'(\mathcal{G})\phi)\dx\\
&=-\int_{\Omega^{*}}(p-1)|g'(\mathcal{G})|^{p-2}g''(\mathcal{G})|\nabla \mathcal{G}|_{\mathcal{A}}^{p}\phi\dx=\int_{\Omega^{*}}Wv^{p-1}\phi\dx,
\end{align*}
where in the penultimate step, we also used~$|g'(\mathcal{G})|^{p-2}g'(\mathcal{G})\phi\in C^{1}_{c}(\Omega^{*})\subseteq W^{1,p}_{c}(\Omega^{*})$ and~$Q'[\mathcal{G}]=0$ in~$\Omega^{*}$.
Consequently,~$v$ is a positive solution of the equation~$Q'_{-W}[u]=0$ in~$\Omega^{*}$. As above, {by the AAP-type theorem (Theorem~\ref{AAP})},~$Q_{-W}$ is nonnegative in~$\Omega^{*}$. In order that~$Q_{-W}$ is critical in~$\Omega^{*}$, by Theorem~\ref{nullc}~(i), it suffices to construct a null-sequence with respect to~$Q_{-W}$ in~$\Omega^{*}$. For all~$k\in\mathbb{N}$ ($k\geq 2$), let
\begin{align*}
\varphi_{k}(t)\triangleq\begin{cases}
				0&0\leq t\leq \frac{1}{k^{2}},\\
			2+\frac{\log t}{\log k}&\frac{1}{k^{2}}\leq t\leq \frac{1}{k},\\
  1&\frac{1}{k}\leq t,
  \end{cases}
\end{align*}
$w_{k}\triangleq\varphi_{k}(v)$, and~$u_{k}\triangleq vw_{k}$. We conclude from \eqref{bcons} that$$\lim_{x\rightarrow0}v(x)=\lim_{x\rightarrow\overline{\infty}}v(x)=0.$$ Therefore, for all~$k\in\mathbb{N}$ ($k\geq 2$),~$w_{k}$ and~$u_{k}$ are compactly supported in~$\Omega^{*}$. It is a simple matter to check that~$\varphi_{k}$ is Lipschitz continuous on~$[0,\infty)$ for all~$k\in\mathbb{N}$ ($k\geq 2$). Since~$v\in C^{1}(\Omega^{*})$,~$w_{k}$ is locally Lipschitz continuous in~$\Omega^{*}$ and hence belongs to~$W^{1,\infty}_{\loc}(\Omega^{*})$ for all~$k\in\mathbb{N}$ ($k\geq 2$). Then~$w_{k}\in W^{1,p}(\Omega^{*})\cap C_{c}(\Omega^{*})$ for all~$k\in\mathbb{N}$ ($k\geq 2$). Since~$v\in C^{1}(\Omega^{*})$, we further deduce that~$u_{k}\in W^{1,p}(\Omega^{*})\cap C_{c}(\Omega^{*})$ for all~$k\in\mathbb{N}$ ($k\geq 2$). Clearly,~$w_{k}$ and~$u_{k}$ are nonnegative in~$\Omega^{*}$ for all~$k\in\mathbb{N}$ ($k\geq 2$). {Note that~$v^{-1}((1/4,3/4))$ is open,~$\emptyset\neq v^{-1}((1/4,3/4))\Subset\Omega^{*}$, and for all~$k\in\mathbb{N}$ ($k\geq 4$),~$\int_{v^{-1}((1/4,\,3/4))}u_{k}^{p}\dx\geq(1/4)^{p}\vol(v^{-1}((1/4,3/4)))>0$.} By Corollary~\ref{co1}, we can find a positive constant~$C$, independent of~$k$, such that for all~$k\in\mathbb{N}$ ($k\geq 2$),
\begin{align}\label{se2}
Q_{-W}[u_{k};\Omega^{*}]\leq C\left(X(v,w_{k})+X(v,w_{k})^{2/p}X(w_{k},v)^{1-2/p}\right),
\end{align}
where~$X(v,w_{k})=\int_{\Omega^{*}}v^{p}|\nabla w_{k}|^{p}_{\mathcal{A}}\dx$ and~$X(w_{k},v)=\int_{\Omega^{*}}w_{k}^{p}|\nabla v|^{p}_{\mathcal{A}}\dx$. By Lemma~\ref{coarea}, there exists a positive constant~$C'$ such that for all sufficiently large~$k\in\mathbb{N}$,
\begin{align*}
&X(v,w_{k})=\int_{\Omega^{*}}v^{p}|\nabla w_{k}|^{p}_{\mathcal{A}}\dx=\int_{\Omega^{*}}v^{p}|\varphi_{k}'(v)|^{p}|\nabla v|^{p}_{\mathcal{A}}\dx\\
&=C'\int_{0}^{1}g(t)^{p}|\varphi_{k}'(g(t))|^{p}|g'(t)|^{p}\dt=C'\int_{[0,1]\cap g^{-1}(1/k^{2},1/k)}g(t)^{p}|\varphi_{k}'(g(t))|^{p}|g'(t)|^{p}\dt\\
&=C'C(p)\left(\frac{1}{\log k}\right)^{p}\int_{[0,1]\cap g^{-1}(1/k^{2},1/k)}\frac{|1-2t|^{p}}{t(1-t)}\dt\\
&\leq C'C(p)\left(\frac{1}{\log k}\right)^{p}\int_{[0,1]\cap g^{-1}(1/k^{2},1/k)}\frac{|1-2t|}{t(1-t)}\dt\\
&= C'C(p)\left(\frac{1}{\log k}\right)^{p-1}
\end{align*}
and \begin{align*}
&X(w_{k},v)=\int_{\Omega^{*}}w_{k}^{p}|\nabla
v|^{p}_{\mathcal{A}}\dx=\int_{\Omega^{*}}\varphi_{k}(v)^{p}|\nabla v|^{p}_{\mathcal{A}}\dx=C'\int_{0}^{1}\varphi_{k}(g(t))^{p}|g'(t)|^{p}\dt\\&=C'C(p)\int_{0}^{1}\varphi_{k}(g(t))^{p}\frac{|1-2t|^{p}}{t(1-t)}\dt\\
&\leq C'C(p)\int_{[0,1]\cap g^{-1}[1/k^{2},(1/4)^{(p-1)/p}]}\varphi_{k}(g(t))^{p}\frac{|1-2t|^{p}}{t(1-t)}\dt\\
&\leq C'C(p)\int_{[0,1]\cap g^{-1}[1/k^{2},(1/4)^{(p-1)/p}]}\frac{|1-2t|}{t(1-t)}\dt\leq C'C(p)\log k,
\end{align*}
where we stipulate that~$\varphi_{k}'(1/k^{2})=\varphi_{k}'(1/k)=0$. By \eqref{se2}, we infer that~$\lim_{k\rightarrow\infty}\\Q_{-W}[u_{k};\Omega^{*}]=0$. {Likewise,~$\{u_{k}/\Vert u_{k}\Vert_{L^{p}(v^{-1}((1/4,\,3/4)))}\}_{k\in\mathbb{N}\,(k\geq 4)}$ is a desired null-sequence.} 
\end{itemize}
By Theorem \ref{nullscon}, up to a positive multiplicative constant,~$v$ is the unique ground state of~$Q_{-W}$ in~$\Omega^{*}$.

Next we demonstrate the null-criticality with respect to~$W$. The demonstration consists of two parts.
\begin{itemize}
\item Suppose that~$1<p\leq n$, or~$p>n$ and~$\sigma=0$. For each~$0<t_{1}<t_{2}<\infty$, by \cite[Theorem~2.32]{ChenCPAM} and checking the proof of Lemma~\ref{coarea}, we may deduce that
\begin{align*}
\int_{\mathcal{G}^{-1}((t_{1},\,t_{2}))}Wv^{p}\dx&=\int_{\mathcal{G}^{-1}((t_{1},\,t_{2}))}\left(\frac{p-1}{p}\right)^{p}\left|\frac{\nabla\mathcal{G}}{\mathcal{G}}\right|_{\mathcal{A}}^{p}\mathcal{G}^{p-1}\dx\\
&=\left(\frac{p-1}{p}\right)^{p}\int_{\mathcal{G}^{-1}((t_{1},\,t_{2}))}\frac{|\nabla\mathcal{G}|_{\mathcal{A}}^{p}}{\mathcal{G}}\dx\\
&=\left(\frac{p-1}{p}\right)^{p}\int_{\mathcal{G}^{-1}((t_{1},\,t_{2}))}\frac{|\nabla\mathcal{G}|_{\mathcal{A}}^{p}}{\mathcal{G}|\nabla\mathcal{G}|}|\nabla\mathcal{G}|\dx\\
&=\left(\frac{p-1}{p}\right)^{p}\int_{t_{1}}^{t_{2}}\int_{\mathcal{G}^{-1}(t)}\frac{|\nabla\mathcal{G}|_{\mathcal{A}}^{p}}{\mathcal{G}|\nabla\mathcal{G}|}\dHnn\dt\\
&=\left(\frac{p-1}{p}\right)^{p}\int_{t_{1}}^{t_{2}}\frac{1}{t}\int_{\mathcal{G}^{-1}(t)}\frac{|\nabla\mathcal{G}|_{\mathcal{A}}^{p}}{|\nabla\mathcal{G}|}\dHnn\dt\\
&=\left(\frac{p-1}{p}\right)^{p}C\int_{t_{1}}^{t_{2}}\frac{1}{t}\dt\\
&\rightarrow\infty\quad\mbox{as}\quad t_{1}\rightarrow 0.
\end{align*}
Then~$\int_{\Omega^{*}}Wv^{p}\dx=\infty$.
\item Suppose that~$p>n$ and~$\sigma>0$. For each~$0<t<\sigma/4$, by \cite[Theorem~2.32]{ChenCPAM} and reviewing the reasoning behind Lemma~\ref{coarea}, we may calculate:
\begin{align*}
&\int_{\mathcal{G}^{-1}((t,\,\sigma/4))}Wv^{p}\dx\\
&=\left(\frac{p-1}{p}\right)^{p}\int_{\mathcal{G}^{-1}((t,\,\sigma/4))}\frac{|\nabla\mathcal{G}|^{p}_{\mathcal{A}}}{|\nabla\mathcal{G}|\mathcal{G}(\sigma-\mathcal{G})}|\sigma-2\mathcal{G}|^{p-2}(2(p-2)\mathcal{G}(\sigma-\mathcal{G})+\sigma^{2})|\nabla\mathcal{G}|\dx\\
&=\left(\frac{p-1}{p}\right)^{p}\int_{t}^{\sigma/4}\int_{\mathcal{G}^{-1}(s)}\frac{|\nabla\mathcal{G}|^{p}_{\mathcal{A}}|\sigma-2\mathcal{G}|^{p-2}}{|\nabla\mathcal{G}|\mathcal{G}(\sigma-\mathcal{G})}(2(p-2)\mathcal{G}(\sigma-\mathcal{G})+\sigma^{2})\dHnn\ds\\
&\geq\frac{\sigma^{p-1}}{2^{p-2}}\left(\frac{p-1}{p}\right)^{p}\int_{t}^{\sigma/4}\frac{1}{s}\int_{\mathcal{G}^{-1}(s)}\frac{|\nabla\mathcal{G}|^{p}_{\mathcal{A}}}{|\nabla\mathcal{G}|}\dHnn\ds\\
&=\frac{\sigma^{p-1}}{2^{p-2}}\left(\frac{p-1}{p}\right)^{p}C\int_{t}^{\sigma/4}\frac{1}{s}\ds\\
&\rightarrow\infty\quad\mbox{as}\quad t\rightarrow 0.
\end{align*}
Therefore,~$\int_{\Omega^{*}}Wv^{p}\dx=\infty$.\qedhere
\end{itemize}
\eproof
Now we present some examples fulfilling the conditions \eqref{bcons}. See \cite[Remark~1.6]{Devyver1} for the case of the~$p$-Laplace equation.
\begin{examples}\label{fexas}
\emph{For simplicity, we only consider~$\hat{x}=0$.
\begin{enumerate}
\item[(i)] Let~$\Omega$ be a bounded domain such that~$0\in\Omega$ and that $\Omega$ satisfies an exterior cone condition or even~$\Omega^{c}$ satisfies a corkscrew condition (see \cite[p.~123]{HKM}). For all~$x\in\Omega^{*}$ and all~$\xi\in\R^{n}$, let~$H(x,\xi)=|x|^{\delta/p}N(\xi)$, where~$\delta\in [0,\infty)$ and $N$ is a uniformly convex norm on~$\R^{n}$ which is differentiable in~$\R^{n}\setminus\{0\}$. Then~$H$ satisfies Assumptions~\ref{ass9} in~$\Omega^{*}$. Suppose, in addition, that~$N$ realizes Assumptions~\ref{C1alpha} and~\ref{assup} in~$\Omega^{*}$. Then~$H$ realizes Assumptions~\ref{C1alpha} and~\ref{assup} in~$\Omega^{*}$. Now let~$p\leq n$. By \cite[Example~2.22 and Theorem~7.39]{HKM}, there exists a Finsler $p$-harmonic function~$\mathcal{G}$ in~$\Omega^{*}$ such that~$\lim_{x\rightarrow 0}\mathcal{G}(x)=\infty$ and~$\lim_{x\rightarrow \overline{\infty}}\mathcal{G}(x)=0$. Checking the proof of \cite[Theorem 7.39]{HKM}, we see that~$\mathcal{G}$ is nonnegative in~$\Omega\setminus\overline{B_{r}(0)}$ and positive on~$S_{r}(0)$ for some~$r>0$. By virtue of the weak comparison principle \cite[Theorem 4.25]{HPR} and the strong maximum principle \cite[6.5]{HKM}, we may show that~$\mathcal{G}$ is positive in~$B_{r}(0)\setminus\{0\}$ and~$\Omega\setminus\overline{B_{r}(0)}$, respectively. Therefore,~$\mathcal{G}$ is positive in~$\Omega\setminus\{0\}$.
\item[(ii)] Suppose that~$0\in\Omega$ and that~$p>n$. Take any Lipschitz domain~$\omega\Subset\Omega$ with~$0\in\omega$. By \cite[Theorem~7.4]{HPR}, we may find some~$\mathcal{G}\in\mathcal{M}_{\omega;\{0\}}^{0}$. By Theorem~\ref{perturbation}~(ii),~$\lambda_{1}(Q_{0};\omega)>0$ (see Definition~\ref{gpeigen}). By Theorem~\ref{peigen}, we may obtain a positive solution~$G\in W^{1,p}_{0}(\omega)$ of~$Q'_{0}[u]=1$ in~$\omega$. By the embedding theory of Sobolev spaces,~$G$ is globally H\"older continuous on~$\overline{\omega}$. Then~$\lim_{x\rightarrow\overline{\infty}}G(x)=0$, where~$\overline{\infty}$ is the ideal point in the one-point compactification of~$\omega$. Take a ball~$B_{r}(0)\Subset\omega$ ($r>0$). Since~$\mathcal{G}\in\mathcal{M}_{\omega;\{0\}}^{0}$, there exists a positive constant~$C$ such that~$\mathcal{G}\leq CG$ in~$\omega\setminus \overline{B_{r}(0)}$. Therefore,~$\lim_{x\rightarrow\overline{\infty}}\mathcal{G}(x)=0$. Moreover, it is impossible that~$\lim_{x\rightarrow 0}\mathcal{G}(x)=0$. Otherwise, 
fix some~$\tau_{0}>0$ such that~$B_{\tau_{0}}(0)\Subset\omega$. For each~$\varepsilon>0$, there exists~$0<\tau<\tau_{0}$ such that~$\mathcal{G}\leq\varepsilon$ on~$S_{\tau}(0)$. It follows that~$\mathcal{G}\leq\varepsilon$ in~$\omega\setminus\overline{B_{\tau}(0)}$ and hence in~$\omega\setminus\overline{B_{\tau_{0}}(0)}$. Letting~$\varepsilon\rightarrow 0$, we see that~$\mathcal{G}\equiv 0$ in~$\omega\setminus\overline{B_{\tau_{0}}(0)}$. But~$\mathcal{G}$ is positive in~$\omega\setminus\{0\}$, a contradiction.
\newline
\newline Now for simplicity, suppose further that the norm family~$H(x,\xi)$ does not depend on~$x$ on~$\omega$. Take~$\tau_{0}>0$ still as above. Clearly,~$\mathcal{G}_{\tau_{0}}(x)\triangleq \mathcal{G}(\tau_{0}x)$ is a positive solution of~$Q'_{0}[u]=0$ in~$B_{1}(0)\setminus\{0\}$ (cf. the proof of \cite[Theorem~9.2~(2)]{Hou}). By \cite[Lemma~2]{Serrin65} (cf. \cite[Chapters~6 and 9]{HKM}),~$\lim_{x\rightarrow 0}\mathcal{G}_{\tau_{0}}(x)=\sigma$ for some~$\sigma\in[0,\infty)$. Thus~$\lim_{x\rightarrow 0}\mathcal{G}(x)=\lim_{x\rightarrow 0}\mathcal{G}_{\tau_{0}}(x/\tau_{0})=\sigma$. Then~$\sigma>0$ because~$\lim_{x\rightarrow 0}\mathcal{G}(x)=0$ is forbidden.
\item[(iii)] Let~$\Omega=\R^{n}$ and suppose that the norm family~$H(x,\cdot)$ does not depend on~$x\in\R^{n}$. Recall that for all~$x\in\R^{n}$,$$H_{0}(x)=\sup_{\xi\in\R^{n}\setminus\{0\}}\frac{x\cdot\xi}{H(\xi)}$$ is the dual norm of~$H$. Then~$\mathcal{G}(x)=H_{0}(x)^{(p-n)/(p-1)}(p\neq n~\mbox{and}~x\neq 0)$ realizes \eqref{bcons} for~$1<p<n$, or~$p>n$ with~$\sigma=0$, respectively. Now we show that~$\mathcal{G}$ is Finsler $p$-harmonic in~$\R^{n}\setminus\{0\}$. Since~$H$ is uniformly convex, by \cite[p.~190]{strictc} and \cite[(3.6)]{Cianchi},~$H_{0}\in C^{1}(\R^{n}\setminus\{0\})$. The ensuing proof is similar to the counterpart of \cite[Lemma~4.2]{Giri1}. See also \cite[Theorem~3.1]{remark}. For each~$x\in\R^{n}\setminus\{0\}$, we compute:
\begin{align*}
\nabla\mathcal{G}(x)=\frac{p-n}{p-1}H_{0}(x)^{(1-n)/(p-1)}\nabla H_{0}(x).
\end{align*}
By \cite[p.~1139]{Jaros}, for all~$x\in\R^{n}\setminus\{0\}$ and all~$t\neq 0$,
\begin{align*}
\nabla H(tx)=\sgn(t)\nabla H(x)\quad\mbox{and}\quad H(\nabla H_{0}(x))=1,
\end{align*}
and for all $x\in\R^{n}$,
\begin{align*}
H\left(H_{0}(x)\nabla H_{0}(x)\right)\nabla H\left(H_{0}(x)\nabla H_{0}(x)\right)=x,
\end{align*} where we define both $H(0)\nabla H(0)$ and $H_{0}(0)\nabla H_{0}(0)$ as $0$.
Then for each~$x\in\R^{n}\setminus\{0\}$,
\begin{align*}
&\mathcal{A}(\nabla\mathcal{G}(x))\\
&=H(\nabla\mathcal{G}(x))^{p-1}\nabla H(\nabla\mathcal{G}(x))\\
&=H\left(\frac{p-n}{p-1}H_{0}(x)^{(1-n)/(p-1)}\nabla H_{0}(x)\right)^{p-1}\nabla H\left(\frac{p-n}{p-1}H_{0}(x)^{(1-n)/(p-1)}\nabla H_{0}(x)\right)\\
&=\left|\frac{p-n}{p-1}\right|^{p-1}H\left(H_{0}(x)^{(1-n)/(p-1)}\nabla H_{0}(x)\right)^{p-2}\times\\
&~~~~H\left(H_{0}(x)^{(1-n)/(p-1)}\nabla H_{0}(x)\right)\nabla H\left(\frac{p-n}{p-1}H_{0}(x)^{(1-n)/(p-1)}\nabla H_{0}(x)\right)\\
&=\left|\frac{p-n}{p-1}\right|^{p-1}H_{0}(x)^{(1-n)(p-2)/(p-1)}H_{0}(x)^{(1-n)/(p-1)-1}\times\\
&~~~~H\left(H_{0}(x)\nabla H_{0}(x)\right)\nabla H\left(\frac{p-n}{p-1}H_{0}(x)^{(1-n)/(p-1)}\nabla H_{0}(x)\right)\\
&=\sgn(p-n)\left|\frac{p-n}{p-1}\right|^{p-1}H_{0}(x)^{-n}H\left(H_{0}(x)\nabla H_{0}(x)\right)\nabla H\left(H_{0}(x)\nabla H_{0}(x)\right)\\
&=\left(\frac{p-n}{p-1}\right)\left|\frac{p-n}{p-1}\right|^{p-2}H_{0}(x)^{-n}x,
\end{align*}
and
\begin{align*}
\dive\mathcal{A}(\nabla\mathcal{G}) &=\left(\frac{p-n}{p-1}\right)\left|\frac{p-n}{p-1}\right|^{p-2}\sum_{i=1}^{n}\frac{\partial (H_{0}(x)^{-n}x_{i})}{\partial x_{i}}\\
&=\left(\frac{p-n}{p-1}\right)\left|\frac{p-n}{p-1}\right|^{p-2}\sum_{i=1}^{n}\left(H_{0}(x)^{-n}-nH_{0}(x)^{-n-1}x_{i}\frac{\partial H_{0}(x)}{\partial x_{i}}\right)\\
&=\left(\frac{p-n}{p-1}\right)\left|\frac{p-n}{p-1}\right|^{p-2}(nH_{0}(x)^{-n}-nH_{0}(x)^{-n-1}x\cdot\nabla H_{0}(x))\\
&=\left(\frac{p-n}{p-1}\right)\left|\frac{p-n}{p-1}\right|^{p-2}(nH_{0}(x)^{-n}-nH_{0}(x)^{-n})=0.
\end{align*}
Hence~$\mathcal{G}$ is Finsler $p$-harmonic in~$\R^{n}\setminus\{0\}$.
\end{enumerate}}
\end{examples}
\section{Optimal Hardy-weights: Nonzero potentials}\label{nzpotential}
In this section, we present optimal Hardy-weights for the functional~$Q$ with nonzero potentials. We also provide two relevant examples of optimal Hardy-weights.

The basis of our construction is Green potentials defined below.
\bdefinition\label{gpotential}
\emph{Suppose that~$Q$ is subcritical in~$\Omega$. For every~$0\leq\phi\in\core\setminus\{0\}$ such that~$\Omega\setminus\supp\phi$ is a domain, a positive solution~$G_{\phi}$ of $Q'[u]=\phi$ in~$\Omega$ such that~$G_{\phi}\in\mathcal{M}_{\Omega;\,\supp\phi}$ is called a \emph{Green potential} of~$Q'$ in~$\Omega$ with density~$\phi$.}
\edefinition
\bremark
\emph{Such a Green potential~$G_{\phi}$ is always continuously differentiable in~$\Omega$ by Remark~\ref{condiffe}.}
\eremark
\bdefinition
\emph{Let~$\omega\Subset\Omega$ be a domain and let~$q$ be as in Definition \ref{Morreydef1}. For $p<n$, the space~$\mathbb{M}^{q}(p;\omega)$ consists of all functions $f$ in~$M^{q}(p;\omega)$ such that~$|f|^{\left(\frac{p^{*}-1}{p-1}\right)'}\in M^{q}(p;\omega)$; for $p=n$, $\mathbb{M}^{q}(p;\omega)$ is defined as $L^{\rho}(\omega)$ for some~$\rho>1$; for $p>n$, $\mathbb{M}^{q}(p;\omega)$ is defined as $M^{q}(p;\omega)=L^{1}(\omega)$. Furthermore, we define $$\mathbb{M}^{q}_{\loc}(p;\Omega)\triangleq \bigcap_{\substack{\omega\Subset\Omega\\\omega~\mbox{is a domain}}}\mathbb{M}^{q}(p;\omega).$$} 
                    \edefinition
                      \bremark
                  \emph{Note that$$\left(\frac{p^{*}-1}{p-1}\right)'=\frac{pn-n+p}{p^{2}}.$$ It is a simple matter to check that~$L^{\rho}(\omega)\subseteq M^{q}(n;\omega)$ for all~$\rho>1$. Then for~$1<p<\infty$,~$\mathbb{M}^{q}(p;\omega)\subseteq M^{q}(p;\omega)$ and~$\mathbb{M}_{\loc}^{q}(p;\Omega)\subseteq M_{\loc}^{q}(p;\Omega)$. Moreover, for~$p<n$,~$L^{q\left(\frac{p^{*}-1}{p-1}\right)'}(\omega)\subseteq\mathbb{M}^{q}(p;\omega)\subseteq L^{\left(\frac{p^{*}-1}{p-1}\right)'}(\omega)$ and~$L_{\loc}^{q\left(\frac{p^{*}-1}{p-1}\right)'}(\Omega)\subseteq\mathbb{M}_{\loc}^{q}(p;\Omega)\subseteq L_{\loc}^{\left(\frac{p^{*}-1}{p-1}\right)'}(\Omega)$.} \eremark
             \bremark\label{nrem}
\emph{Let $\omega\Subset\Omega$ be a Lipschitz domain. Let $H$ satisfy Assumptions \ref{ass9} and let $V\in\mathbb{M}^{q}(p;\omega)$.  Let $g\in M^{q}(p;\omega)$ if $p\neq n$ and let $g\in L^{\bar{\rho}}(\omega)$ ($\bar{\rho}>1$) if $p=n$. By \cite[Theorem~2]{Liang} for $p<n$ and \cite[Chapter~4, Theorem~7.1]{LAQE} for $p=n$, solutions of the Dirichlet problem\begin{align}\label{bvp}
\begin{cases}
		Q'[u]=g&\mbox{in}~\omega,\\	u=f&\mbox{on}~\partial\omega,
		\end{cases}
\end{align} are bounded in~$\omega$, where $f\in W^{1,p}(\omega)\cap C(\overline{\omega})$. By \cite[Theorem~4.11 and Corollary~4.18]{Maly97}, solutions of \eqref{bvp} are continuous on $\overline{\omega}$. For $p>n$, by the embedding theory of Sobolev spaces, solutions of \eqref{bvp} are also continuous on $\overline{\omega}$.}

\emph{Note that in \cite[Theorems 4.24 and 4.25]{HPR}, it should have been further assumed that~$V\in\mathbb{M}^{q}(p;\omega)$ and if $p=n$, $g\in L^{\bar{\rho}}(\omega)$ ($\bar{\rho}>1$), while in \cite[Theorems 7.2, 7.4, 7.7, 7.8 and 7.9]{HPR}, we need~$V\in\mathbb{M}_{\loc}^{q}(p;\Omega)$ (instead of~$V\in\widetilde{M}_{\loc}^{q}(p;\Omega)$ in \cite[arXiv]{HPR}).}
\eremark
                    Now we show the existence of Green potentials when $V\in\mathbb{M}^{q}_{\loc}(p;\Omega)\cap\widehat{M}^{q}_{\loc}(p;\Omega)$.
\blemma\label{gpcon}
Let~$V\in\mathbb{M}^{q}_{\loc}(p;\Omega)\cap\widehat{M}^{q}_{\loc}(p;\Omega)$. Suppose that~$Q$ is subcritical in~$\Omega$. For every~$0\leq\phi\in\core\setminus\{0\}$ such that~$\Omega\setminus\supp\phi$ is a domain, there exists a Green potential~$G_{\phi}$ of~$Q'$ in~$\Omega$ with density~$\phi$.
\elemma
\bproof
The proof is based on an exhaustion argument, similar to those of \cite[Lemma~3.2]{Versano} and \cite[Theorem~9.7]{Hou}. Let~$\{\omega_{k}\}_{k\in\mathbb{N}}$ be a smooth exhaustion of~$\Omega$ such that~$\supp\phi\subseteq\omega_{1}$. By Theorem~\ref{perturbation}~(ii),~$\lambda_{1}(Q;\omega_{k})>0$ (see Definition~\ref{gpeigen}) for all~$k\in\mathbb{N}$. By Theorem~\ref{peigen}, for every~$k\in\mathbb{N}$, there exists a positive solution~$G_{k}\in W^{1,p}_{0}(\omega_{k})$ of~$Q'[u]=\phi$ in~$\omega_{k}$. By Remark~\ref{nrem},~$G_{k}$ is continuous on~$\overline{\omega_{k}}$ for all~$k\in\mathbb{N}$. By the weak comparison principle \cite[Theorem~4.25]{HPR}, the sequence~$\{G_{k}\}_{k\in\mathbb{N}}$ is increasing in~$\Omega$. 

Now we exploit the method of \cite[Theorem~9.7]{Hou} to complete the proof.

If there exists~$x_{0}\in\Omega$ such that~$\lim_{k\rightarrow\infty}G_{k}(x_{0})=\infty$, then let~$\Gamma_{k}\triangleq G_{k}/G_{k}(x_{0})$ for all sufficiently large~$k\in\mathbb{N}$ such that~$x_{0}\in\omega_{k}$. By Harnack's inequality \cite[Theorem~3.14]{Maly97} for~$p\leq n$ and \cite[Theorem~7.4.1]{Serrin} for~$p>n$, for all sufficiently large~$k\in\mathbb{N}$,~$\Gamma_{k}$ is locally uniformly bounded in~$\Omega$. By the Harnack convergence principle \cite[Theorem~3.5]{HPR}, up to a subsequence, in~$\Omega$,~$\Gamma_{k}$ converges locally uniformly to a positive solution~$\Gamma$ of~$Q'[u]=0$ as~$k\rightarrow\infty$. Since $Q$ is subcritical in~$\Omega$, by Theorem~\ref{nullc}~(ii), there exists a positive~$W\in C(\Omega)$ such that~$Q_{V-W}$ is nonnegative in~$\Omega$. By the AAP-type theorem (Theorem~\ref{AAP}), there exists a positive solution~$v$ of~$Q'_{V-W}[u]=0$ in~$\Omega$. Let~$K$ be an arbitrary admissible compact subset of~$\Omega$. Then there exists a positive constant~$c_{K}$ such that~$\Gamma\leq c_{K} v$ on~$K$. For each~$\varepsilon>0$, up to a subsequence, for all sufficiently large~$k\in\mathbb{N}$,
\begin{align*}
\begin{cases}
		Q'[\Gamma_{k}]=\frac{\phi}{G_{k}(x_{0})^{p-1}}\leq ((1+\varepsilon)c_{K}v)^{p-1}W=Q'[(1+\varepsilon)c_{K} v]&\mbox{in}~\omega_{k}\setminus K,\\
			\Gamma_{k}\leq (1+\varepsilon)c_{K}v&\mbox{on}~\partial(\omega_{k}\setminus K).
		\end{cases}
\end{align*}
By Theorem~\ref{perturbation}~(ii),~$\lambda_{1}(Q;\omega_{k}\setminus K)>0$ (see Definition~\ref{gpeigen}) for all sufficiently large~$k\in\mathbb{N}$. By the weak comparison principle \cite[Theorem~4.25]{HPR}, for each~$\varepsilon>0$, up to a subsequence, for all sufficiently large~$k\in\mathbb{N}$,~$\Gamma_{k}\leq (1+\varepsilon)c_{K}v$ in~$\omega_{k}\setminus K$. Therefore,~$\Gamma\leq (1+\varepsilon)c_{K}v$ in~$\Omega\setminus K$. Letting~$\varepsilon\rightarrow 0$, we get~$\Gamma\leq c_{K}v$ in~$\Omega\setminus K$. Then~$\Gamma\leq c_{K}v$ in~$\Omega$. Let~$c_{0}\triangleq\inf\{c>0\,|\,\Gamma\leq cv~\mbox{in}~\Omega\}$. Clearly,~$\Gamma\leq c_{0}v$ in~$\Omega$ and~$c_{0}>0$. Inasmuch as~$\Gamma$ and~$c_{0}v$ are respectively positive solutions of~$Q'[u]=0$ and~$Q'_{V-W}[u]=0$ in~$\Omega$, it is easy to see that~$\Gamma\neq c_{0}v$. Consequently, there exists~$x_{1}\in\Omega$ such that~$\Gamma(x_{1})<c_{0}v(x_{1})$. Furthermore, there exist~$\delta,r>0$ such that~$\Gamma\leq (1-\delta)c_{0}v$ on~$\overline{B_{r}(x_{1})}$. By our previous argument,~$\Gamma\leq (1-\delta)c_{0}v$ in~$\Omega$, which is contradictory to the definition of~$c_{0}$.

Therefore, for all~$x\in\Omega$,~$\{G_{k}(x)\}_{k\in\mathbb{N}}$ is bounded. By Harnack's inequality \cite[Theorem~3.14]{Maly97} for~$p\leq n$ and \cite[Theorem~7.4.1]{Serrin} for~$p>n$,~$\{G_{k}\}_{k\in\mathbb{N}}$ is locally uniformly bounded in~$\Omega$. Then we may use Harnack's convergence principle \cite[Theorem~3.5]{HPR} to deduce that up to a subsequence, in~$\Omega$,~$\{G_{k}\}_{k\in\mathbb{N}}$ converges locally uniformly to a positive solution~$G_{\phi}$ of~$Q'[u]=\phi$. 

Next we show that~$G_{\phi}\in\mathcal{M}_{\Omega;\,\supp\phi}$. Let~$K$ be an arbitrary admissible compact subset of~$\Omega$ with~$\supp\phi\Subset \mathring{K}$ and consider an arbitrary positive solution~$v\in C\left(\Omega\setminus \mathring{K}\right)$  of~$Q'[w]=g$  in~$\Omega\setminus K$ such that $G_{\phi}\leq v$ on~$\partial K$, where~$g\in M^{q}_{\loc}(p;\Omega)$ for $p\neq n$, $g\in L^{\bar{\rho}}_{\loc}(\Omega))$ ($\bar{\rho}>1$) for $p=n$, and $g$ is nonnegative a.e. in $\Omega\setminus K$. Up to a subsequence, for all sufficiently large~$k\in\mathbb{N}$,
\begin{align*}
\begin{cases}
		Q'[G_{k}]=0\leq g=Q'[v]&\mbox{in}~\omega_{k}\setminus K,\\
			G_{k}\leq v&\mbox{on}~\partial(\omega_{k}\setminus K).
		\end{cases}
\end{align*}
By Theorem~\ref{perturbation}~(ii),~$\lambda_{1}(Q;\omega_{k}\setminus K)>0$ (see Definition~\ref{gpeigen}) for all sufficiently large~$k\in\mathbb{N}$. By the weak comparison principle \cite[Theorem~4.25]{HPR}, up to a subsequence, for all sufficiently large~$k\in\mathbb{N}$,~$G_{k}\leq v$ in~$\omega_{k}\setminus K$. It follows that~$G_{\phi}\leq v$ in~$\Omega\setminus K$. 
\eproof
For the construction of optimal Hardy-weights in the case of nonzero potentials, we need the following preparatory result.
\blemma\label{re1}
Suppose that~$Q$ is subcritical in~$\Omega$. Let~$G_{\phi}$ be a Green potential of~$Q'$ in~$\Omega$ with density~$0\leq\phi\in\core\setminus\{0\}$ such that~$\Omega\setminus\supp\phi$ is a domain. Assume the following conditions:
$$\lim_{x\rightarrow\overline{\infty}}G_{\phi}=0,\quad\int_{\Omega}|V|G_{\phi}^{p-1}\dx<\infty,\quad\mbox{and}\quad\int_{\Omega}VG_{\phi}^{p-1}\dx<0~\mbox{or}~V\leq 0~\mbox{in}~\Omega.$$
Then there exist two positive constants~$C_{0}$ and~$M_{\phi}<\sup_{\Omega}G_{\phi}$ such that for almost all~$0<t<\sup_{\Omega}G_{\phi}$,
$$\int_{G_{\phi}^{-1}(t)}\frac{|\nabla G_{\phi}|_{\mathcal{A}}^{p}}{|\nabla G_{\phi}|}\dHnn\leq C_{0},$$ and for almost all~$0<t<M_{\phi}$ with~$\supp\phi\subseteq G_{\phi}^{-1}((t,\infty))$,
$$\frac{1}{C_{0}}\leq\int_{G_{\phi}^{-1}(t)}\frac{|\nabla G_{\phi}|_{\mathcal{A}}^{p}}{|\nabla G_{\phi}|}\dHnn.$$
\elemma
\bproof
The proof is similar to that of \cite[Lemma~3.7]{Versano}. By our assumption,~$\lim_{x\rightarrow\overline{\infty}}G_{\phi}=0$. Hence, for every~$0<t<\sup_{\Omega}G_{\phi}$,~$G_{\phi}^{-1}((t,\infty))\Subset\omega_{t}$ for some domain~$\omega_{t}\Subset\Omega$ and~$G_{\phi}<t$ in~$\Omega\setminus\omega_{t}$. By \cite[Theorem~5.9~(i)]{evansf},~$G_{\phi}^{-1}((t,\infty))$ has finite perimeter in~$\omega_{t}$ for almost all~$0<t<\sup_{\Omega}G_{\phi}$. By the generalized Sard theorem \cite[Theorem~1.2]{Sard}, for almost all~$0<t<\sup_{\Omega}G_{\phi}$,$$\mathcal{H}^{n-1}\left(\{x\in\Omega~|~G_{\phi}(x)=t~\mbox{and}~\nabla G_{\phi}(x)=0\}\right)=0.$$ It follows that~$\partial G_{\phi}^{-1}((t,\infty))=G_{\phi}^{-1}(t)$ for almost all~$0<t<\sup_{\Omega}G_{\phi}$. Then~$G_{\phi}^{-1}((t,\infty))$ is an open set with almost~$C^{1}$-boundary (see \cite[Section~9.3]{Maggi}) for almost all~$0<t<\sup_{\Omega}G_{\phi}$. 
Note that~$G_{\phi}^{-1}((t,\infty))=\{x\in\omega_{t}~|~G_{\phi}>t\}$. Since~$-\dive\mathcal{A}(x,\nabla G_{\phi})+VG_{\phi}^{p-1}=\phi$ in~$\Omega$, by \cite[Theorems~5.2 and 7.2]{ChenCPAM}, we may deduce that for almost all~$0<t<\sup_{\Omega}G_{\phi}$, 
\begin{align}\label{eq1}
&\int_{G_{\phi}^{-1}((t,\infty))}\left(\phi-VG_{\phi}^{p-1}\right)\dx=-\int_{\partial G_{\phi}^{-1}((t,\infty))}\mathcal{A}(x,\nabla G_{\phi})\cdot\frac{-\nabla G_{\phi}}{|\nabla G_{\phi}|}\dHnn\notag\\
&=\int_{G_{\phi}^{-1}(t)}\frac{|\nabla G_{\phi}|_{\mathcal{A}}^{p}}{|\nabla G_{\phi}|}\dHnn.
\end{align}
Take a positive constant~$C_{0}\geq\max\{\int_{\Omega}\left(\phi+|V|G_{\phi}^{p-1}\right)\dx,1/\int_{\Omega}\phi\dx\}$.
Then for all~$0<t<\sup_{\Omega}G_{\phi}$,\begin{align}\label{eq2}
\int_{G_{\phi}^{-1}((t,\infty))}(\phi-VG_{\phi}^{p-1})\dx\leq\int_{\Omega}\left(\phi+|V|G_{\phi}^{p-1}\right)\dx\leq C_{0}.
\end{align}
Since~$\int_{\Omega}VG_{\phi}^{p-1}\dx<0$ or~$V\leq 0$ in~$\Omega$, there exists~$0<M_{\phi}<\sup_{\Omega}G_{\phi}$ such that for all~$0<t<M_{\phi}$,$$\int_{G_{\phi}^{-1}((t,\infty))}VG_{\phi}^{p-1}\dx\leq 0,$$ which is proved by contradiction when~$\int_{\Omega}VG_{\phi}^{p-1}\dx<0$. Therefore, for all~$0<t<M_{\phi}$ with~$\supp\phi\subseteq G_{\phi}^{-1}((t,\infty))$,\begin{align}\label{eq3}
\frac{1}{C_{0}}\leq\int_{\Omega}\phi\dx=\int_{G_{\phi}^{-1}((t,\infty))}\phi\dx\leq\int_{G_{\phi}^{-1}((t,\infty))}(\phi-VG_{\phi}^{p-1})\dx.\end{align} 
The desired conclusions follow from \eqref{eq1}, \eqref{eq2}, and \eqref{eq3}.
\eproof
Our construction of optimal Hardy-weights is done in the following theorem. Recall that~$c_{p}=(p/(p-1))^{p-1}$ and that~$f_{0}(t)=t^{(p-1)/p}$ for all $t>0$.
\btheorem\label{mainoh}
Suppose that~$Q_{c_{p}V}$ is subcritical in~$\Omega$. Let~$G_{\phi}$ be a Green potential of~$Q'_{c_{p}V}$ in~$\Omega$ with density~$0\leq\phi\in\core\setminus\{0\}$ such that~$\Omega\setminus\supp\phi$ is a domain. Suppose that Assumption~\ref{assup} and the following conditions hold:$$\lim_{x\rightarrow\overline{\infty}}G_{\phi}=0,\quad\int_{\Omega}|V|G_{\phi}^{p-1}\dx<\infty,\quad\mbox{and}\quad\int_{\Omega}VG_{\phi}^{p-1}\dx<0~\mbox{or}~V\leq 0~\mbox{in}~\Omega.$$
Let\begin{align*}
    W&\triangleq\frac{-(p-1)f'_{0}(G_{\phi})^{p-2}f_{0}''(G_{\phi})|\nabla G_{\phi}|_{\mathcal{A}}^{p}+f'_{0}(G_{\phi})^{p-1}\phi}{f_{0}(G_{\phi})^{p-1}}\\
    &=\left(\frac{p-1}{p}\right)^{p}\frac{|\nabla G_{\phi}|_{\mathcal{A}}^{p}}{G_{\phi}^{p}}+\left(\frac{p-1}{p}\right)^{p-1}G_{\phi}^{1-p}\phi.
\end{align*}Then~$0\leq W\in L^{\infty}_{\loc}(\Omega)\subseteq\widehat{M}^{q}_{\loc}(p;\Omega)$ and~$Q_{V-W}$ is critical in~$\Omega$ with a ground state~$f_{0}(G_{\phi})$ (up to a positive multiplicative constant) such that~$\int_{\Omega}Wf_{0}(G_{\phi})^{p}\dx=\infty$. In particular,~$W$ is an optimal Hardy-weight of~$Q$ in~$\Omega$.
\etheorem
\bproof
The proof is similar to that of \cite[Lemma~3.9]{Versano}. Because~$f_{0},f_{0}'>0$,~$\phi\geq 0$ in~$\Omega$, and~$f_{0}''<0$, we see at once that~$W$ is nonnegative in~$\Omega$. Recalling that~$G_{\phi}\in C^{1}(\Omega)$ and~$G_{\phi}>0$ in~$\Omega$, we may easily check that~$W\in L^{\infty}_{\loc}(\Omega)$. Since~$Q_{c_{p}V}$ is subcritical in~$\Omega$,~$Q=Q_{V}$ is subcritical in~$\Omega$ by Theorem~\ref{perturbation}~(iii). By Lemma~\ref{re1}, there exist two positive constants~$C_{0}$ and~$M_{\phi}<\sup_{\Omega}G_{\phi}$ such that for almost all~$0<t<\sup_{\Omega}G_{\phi}$,
\begin{align}\label{ubgphi}
\int_{G_{\phi}^{-1}(t)}\frac{|\nabla G_{\phi}|_{\mathcal{A}}^{p}}{|\nabla G_{\phi}|}\dHnn\leq C_{0},
\end{align}
and for almost all~$0<t<M_{\phi}$ with~$\supp\phi\subseteq G_{\phi}^{-1}((t,\infty))$,
\begin{align}\label{lbgphi}
\frac{1}{C_{0}}\leq\int_{G_{\phi}^{-1}(t)}\frac{|\nabla G_{\phi}|_{\mathcal{A}}^{p}}{|\nabla G_{\phi}|}\dHnn.    
\end{align} Since~$\lim_{x\rightarrow\overline{\infty}}G_{\phi}=0$, we may take a smooth domain~$\omega\Subset\Omega$ such that~$\supp\phi\subseteq\omega$ and~$\sup_{\Omega\setminus\overline{\omega}}G_{\phi}<\min\{1,M_{\phi}\}$.

As in the proof of Theorem~\ref{zerop}, for every~$k\in\mathbb{N}$ ($k\geq 2$), let
\begin{align*}
\varphi_{k}(t)\triangleq\begin{cases}
				0&0\leq t\leq \frac{1}{k^{2}},\\
			2+\frac{\log t}{\log k}&\frac{1}{k^{2}}\leq t\leq \frac{1}{k},\\
  1&\frac{1}{k}\leq t\leq k,\\
  2-\frac{\log t}{\log k}&k\leq t\leq k^{2},\\
  0& t\geq k^{2},
		\end{cases}
\end{align*}
$w_{k}\triangleq\varphi_{k}(f_{0}(G_{\phi}))$,
and~$u_{k}\triangleq w_{k}f_{0}(G_{\phi})$. Corollary~\ref{co2} implies that~$f_{0}(G_{\phi})$ is a positive solution of~$Q'_{V-W}[u]=0$ in~$\Omega$. The AAP-type theorem (Theorem~\ref{AAP}) further ensures that~$Q_{V-W}$ is nonnegative in~$\Omega$. Observe that~$\varphi_{k}$ is Lipschitz continuous on~$[0,\infty)$ for all~$k\in\mathbb{N}$ ($k\geq 2$). It follows that~$w_{k}$ is locally Lipschitz continuous in~$\Omega$ and thus belongs to~$W^{1,\infty}_{\loc}(\Omega)$ for all~$k\in\mathbb{N}$ ($k\geq 2$). Because~$\lim_{x\rightarrow\overline{\infty}}G_{\phi}=0$,~$w_{k}$ has compact support in~$\Omega$ for all~$k\in\mathbb{N}$ ($k\geq 2$). Consequently,~$w_{k}\in W^{1,p}(\Omega)\cap C_{c}(\Omega)$ for all~$k\in\mathbb{N}$ ($k\geq 2$). Then~$u_{k}$ also belongs to~$W^{1,p}(\Omega)\cap C_{c}(\Omega)$ for all~$k\in\mathbb{N}$ ($k\geq 2$) because~$f_{0}(G_{\phi})\in C^{1}(\Omega)$. Both~$w_{k}$ and~$u_{k}$ are nonnegative in~$\Omega$ for all~$k\in\mathbb{N}$ ($k\geq 2$). By Corollary~\ref{co1}, there exists a positive constant~$C$, independent of~$k$, such that for all~$k\in\mathbb{N}$ ($k\geq 2$),
\begin{align*}
Q_{V-W}[u_{k}]\leq \begin{cases}		CX(f_{0}(G_{\phi}),w_{k})&~\mbox{if}~p\leq 2,\\
		C\left(X(f_{0}(G_{\phi}),w_{k})+X(f_{0}(G_{\phi}),w_{k})^{2/p}X(w_{k},f_{0}(G_{\phi}))^{1-2/p}\right)&~\mbox{if}~p>2,
		\end{cases}
\end{align*}
where~$X(f_{0}(G_{\phi}),w_{k})\triangleq\int_{\Omega}f_{0}(G_{\phi})^{p}|\nabla w_{k}|^{p}_{\mathcal{A}}\dx$ and~$X(w_{k},f_{0}(G_{\phi}))\triangleq\int_{\Omega}w_{k}^{p}|\nabla (f_{0}(G_{\phi}))|^{p}_{\mathcal{A}}\dx$.
There exists a positive integer~$N$ such that for all~$k>N$, $1/k\leq f_{0}(G_{\phi})\leq k$ in some domain~$\omega'\Subset\Omega$ with~$\omega\Subset\omega'$. By enlarging~$N$, we may ensure~$f_{0}(\sup_{\omega'}G_{\phi})\leq N$. By the definition of~$\{\varphi_{k}\}_{k\in\mathbb{N}}$, for all~$k>N$,~$\nabla w_{k}=0$ on~$\overline{\omega}$.
Thus for all~$k>N$,
\begin{align*}
X(f_{0}(G_{\phi}),w_{k})=\int_{\Omega\setminus\overline{\omega}}f_{0}(G_{\phi})^{p}|\nabla w_{k}|^{p}_{\mathcal{A}}\dx.
\end{align*}
For every~$k\in\mathbb{N}$ ($k\geq 2$), suppose that~$\supp w_{k}\subseteq\omega_{k}$ for some smooth domain~$\omega_{k}\Subset\Omega$. We may also assume that~$\omega\Subset\omega_{k}$ for all~$k\in\mathbb{N}$ ($k\geq 2$). Then by \cite[Theorem~2.32]{ChenCPAM}, for all~$k>N$,
\begin{align*}
X(f_{0}(G_{\phi}),w_{k})&=\int_{\omega_{k}\setminus\overline{\omega}}f_{0}(G_{\phi})^{p}|\nabla w_{k}|^{p}_{\mathcal{A}}\dx\\
&=\int_{\omega_{k}\setminus\overline{\omega}}f_{0}(G_{\phi})^{p}|\varphi_{k}'(f_{0}(G_{\phi}))|^{p}f_{0}'(G_{\phi})^{p}\frac{|\nabla G_{\phi}|^{p}_{\mathcal{A}}}{|\nabla G_{\phi}|}|\nabla G_{\phi}|\dx\\
&=\int_{\inf_{\omega_{k}\setminus\overline{\omega}}G_{\phi}}^{\sup_{\omega_{k}\setminus\overline{\omega}}G_{\phi}}f_{0}(t)^{p}|\varphi_{k}'(f_{0}(t))|^{p}f_{0}'(t)^{p}\dt\int_{\{x\in\omega_{k}\setminus\overline{\omega}\,|\,G_{\phi}(x)=t\}}\frac{|\nabla G_{\phi}|^{p}_{\mathcal{A}}}{|\nabla G_{\phi}|}\dHnn,
\end{align*}
where it is understood that~$\varphi_{k}'(1/k^{2})=\varphi_{k}'(1/k)=0$. 
By \eqref{ubgphi}, for almost all~$\inf_{\omega_{k}\setminus\overline{\omega}}G_{\phi}<t<\sup_{\omega_{k}\setminus\overline{\omega}}G_{\phi}$,
$$\int_{\{x\in\omega_{k}\setminus\overline{\omega}\,|\,G_{\phi}(x)=t\}}\frac{|\nabla G_{\phi}|^{p}_{\mathcal{A}}}{|\nabla G_{\phi}|}\dHnn\leq\int_{G_{\phi}^{-1}(t)}\frac{|\nabla G_{\phi}|_{\mathcal{A}}^{p}}{|\nabla G_{\phi}|}\dHnn\leq C_{0}.$$
We now calculate: for all~$k>N$,
\begin{align*}
&\int_{\inf_{\omega_{k}\setminus\overline{\omega}}G_{\phi}}^{\sup_{\omega_{k}\setminus\overline{\omega}}G_{\phi}}f_{0}(t)^{p}|\varphi_{k}'(f_{0}(t))|^{p}f_{0}'(t)^{p}\dt\leq\int_{(1/k^{2})^{p/(p-1)}}^{(1/k)^{p/(p-1)}}f_{0}(t)^{p}|\varphi_{k}'(f_{0}(t))|^{p}f_{0}'(t)^{p}\dt\\
&=C(p)\int_{1/k^{2}}^{1/k}s^{p-1}|\varphi_{k}'(s)|^{p}\ds=\frac{C(p)}{\log^{p}k}\int_{1/k^{2}}^{1/k}\frac{1}{s}\ds=\frac{C(p)}{\log^{p-1}k}.
\end{align*}
Then for all~$k>N$,$$X(f_{0}(G_{\phi}),w_{k})\leq\frac{C(p)C_{0}}{\log^{p-1}k}.$$ Combining~$\sup_{\Omega\setminus\overline{\omega}}G_{\phi}<1$ with~$f_{0}(\sup_{\omega'}G_{\phi})\leq N$ yields~$\sup_{\omega_{k}}G_{\phi}\leq N^{p/(p-1)}$. We proceed to estimate:
\begin{align*}
&X(w_{k},f_{0}(G_{\phi}))=\int_{\omega_{k}}w_{k}^{p}|\nabla (f_{0}(G_{\phi}))|^{p}_{\mathcal{A}}\dx=\int_{\omega_{k}}\varphi_{k}(f_{0}(G_{\phi}))^{p}f_{0}'(G_{\phi})^{p}|\nabla G_{\phi}|^{p}_{\mathcal{A}}\dx\\
&=\int_{\omega_{k}}\varphi_{k}(f_{0}(G_{\phi}))^{p}f_{0}'(G_{\phi})^{p}\frac{|\nabla G_{\phi}|^{p}_{\mathcal{A}}}{|\nabla G_{\phi}|}|\nabla G_{\phi}|\dx\\
&=\int_{\inf_{\omega_{k}}G_{\phi}}^{\sup_{\omega_{k}}G_{\phi}}\varphi_{k}(f_{0}(t))^{p}f_{0}'(t)^{p}\dt\int_{\{x\in\omega_{k}\,|\,G_{\phi}(x)=t\}}\frac{|\nabla G_{\phi}|^{p}_{\mathcal{A}}}{|\nabla G_{\phi}|}\dHnn\\
&\leq C_{0}\int_{\inf_{\omega_{k}}G_{\phi}}^{\sup_{\omega_{k}}G_{\phi}}\varphi_{k}(f_{0}(t))^{p}f_{0}'(t)^{p}\dt\\
&\leq C_{0}\int_{(1/k^{2})^{p/(p-1)}}^{(1/k)^{p/(p-1)}}\varphi_{k}(f_{0}(t))^{p}f_{0}'(t)^{p}\dt+C_{0}\int_{(1/k)^{p/(p-1)}}^{N^{p/(p-1)}}\varphi_{k}(f_{0}(t))^{p}f_{0}'(t)^{p}\dt\\
&=C_{0}C(p)\int_{(1/k^{2})^{p/(p-1)}}^{(1/k)^{p/(p-1)}}\left(2+\frac{\log f_{0}(t)}{\log k}\right)^{p}\frac{1}{t}\dt+C_{0}C(p)\int_{(1/k)^{p/(p-1)}}^{N^{p/(p-1)}}\frac{1}{t}\dt\\
&\leq C_{0}C(p)\int_{(1/k^{2})^{p/(p-1)}}^{(1/k)^{p/(p-1)}}\frac{1}{t}\dt+C_{0}C(p)\log k\leq C_{0}C(p)\log k.
\end{align*}
Therefore,~$\lim_{k\rightarrow\infty}Q_{V-W}[u_{k}]=0.$ Take~$0<\varepsilon<\min\{\inf_{\overline{\omega}}f_{0}(G_{\phi}),1\}$. Then the preimage $f_{0}(G_{\phi})^{-1}((\varepsilon/2,\varepsilon))\Subset\Omega$ and is nonempty and open. Furthermore, there exists~$N'\in\mathbb{N}$ such that for all~$k\geq N'$,
\begin{align*}
&\int_{f_{0}(G_{\phi})^{-1}((\varepsilon/2,\,\varepsilon))}|u_{k}|^{p}\dx=\int_{f_{0}(G_{\phi})^{-1}((\varepsilon/2,\,\varepsilon))}|\varphi_{k}(f_{0}(G_{\phi}))|^{p}f_{0}(G_{\phi})^{p}\dx\\
&=\int_{f_{0}(G_{\phi})^{-1}((\varepsilon/2,\,\varepsilon))}f_{0}(G_{\phi})^{p}\dx>0.
\end{align*}Thus~$\{u_{k}/\Vert u_{k}\Vert_{L^{p}(f_{0}(G_{\phi})^{-1}((\varepsilon/2,\,\varepsilon)))}\}_{k\geq N'}$ is a null-sequence of~$Q_{V-W}$. By Theorem~\ref{nullc}~(i), $Q_{V-W}$ is critical in~$\Omega$. Since~$f_{0}(G_{\phi})$ is a positive solution of~$Q'_{V-W}[u]=0$ in~$\Omega$, by Theorem \ref{nullscon},~$f_{0}(G_{\phi})$ is a ground state of~$Q_{V-W}$ up to a positive multiplicative constant. 

Note that in~$\Omega\setminus\supp\phi$,
$$W=\left(\frac{p-1}{p}\right)^{p}\frac{|\nabla G_{\phi}|_{\mathcal{A}}^{p}}{G_{\phi}^{p}}.$$ Let~$0<\tau<\inf_{\overline{\omega}}G_{\phi}$. Note that~$G_{\phi}^{-1}((\tau,\inf_{\overline{\omega}}G_{\phi}))\subseteq\Omega\setminus\supp\phi$ because~$\supp\phi\subseteq\omega$ and~$G_{\phi}^{-1}((\tau,\inf_{\overline{\omega}}G_{\phi}))\subseteq\Omega\setminus\overline{\omega}$. Since~$\lim_{x\rightarrow\overline{\infty}}G_{\phi}=0$ and~$\tau>0$, we get~$G_{\phi}^{-1}((\tau,\inf_{\overline{\omega}}G_{\phi}))\Subset\Omega$. Hence~$G_{\phi}$ is Lipschitz in~$G_{\phi}^{-1}((\tau,\inf_{\overline{\omega}}G_{\phi}))$. Then by \cite[Theorem~2.32]{ChenCPAM} and \eqref{lbgphi},
\begin{align*}
&\int_{G_{\phi}^{-1}((\tau,\,\inf_{\overline{\omega}}G_{\phi}))}Wf_{0}(G_{\phi})^{p}\dx=\left(\frac{p-1}{p}\right)^{p}\int_{G_{\phi}^{-1}((\tau,\,\inf_{\overline{\omega}}G_{\phi}))}\frac{|\nabla G_{\phi}|_{\mathcal{A}}^{p}}{G_{\phi}^{p}|\nabla G_{\phi}|}f_{0}(G_{\phi})^{p}|\nabla G_{\phi}|\dx\\
&=\left(\frac{p-1}{p}\right)^{p}\int_{\tau}^{\inf_{\overline{\omega}}G_{\phi}}\frac{{f_{0}(t)}^{p}}{t^{p}}\dt\int_{\left\{x\in G_{\phi}^{-1}((\tau,\,\inf_{\overline{\omega}}G_{\phi}))\,|\,G_{\phi}(x)=t\right\}}\frac{|\nabla G_{\phi}|_{\mathcal{A}}^{p}}{|\nabla G_{\phi}|}\dHnn\\
&=\left(\frac{p-1}{p}\right)^{p}\int_{\tau}^{\inf_{\overline{\omega}}G_{\phi}}\frac{{f_{0}(t)}^{p}}{t^{p}}\dt\int_{G^{-1}_{\phi}(t)}\frac{|\nabla G_{\phi}|_{\mathcal{A}}^{p}}{|\nabla G_{\phi}|}\dHnn\\
&\geq\frac{1}{C_{0}}\left(\frac{p-1}{p}\right)^{p}\int_{\tau}^{\inf_{\overline{\omega}}G_{\phi}}\frac{{f_{0}(t)}^{p}}{t^{p}}\dt=\frac{1}{C_{0}}\left(\frac{p-1}{p}\right)^{p}\int_{\tau}^{\inf_{\overline{\omega}}G_{\phi}}\frac{1}{t}\dt.
\end{align*}
With sufficiently small~$\tau$, we see that$$\int_{G_{\phi}^{-1}((0,\,\inf_{\overline{\omega}}G_{\phi}))}Wf_{0}(G_{\phi})^{p}\dx=\infty.$$
Immediately, we get~$\int_{\Omega}Wf_{0}(G_{\phi})^{p}\dx=\infty$.
\eproof
With Theorem~\ref{mainoh} at our disposal, the proof of the following corollary is trivial and hence omitted. See also \cite[Theorem~1.1]{Versano}.
\bcorollary\label{reph}
Suppose that~$Q$ is subcritical in~$\Omega$. Let~$G_{\phi}$ be a Green potential of~$Q'$ in~$\Omega$ with density~$0\leq\phi\in\core\setminus\{0\}$ such that~$\Omega\setminus\supp\phi$ is a domain. Suppose that Assumption~\ref{assup} and the following conditions hold: $$\lim_{x\rightarrow\overline{\infty}}G_{\phi}=0,\quad\int_{\Omega}|V|G_{\phi}^{p-1}\dx<\infty,\quad\mbox{and}\quad\int_{\Omega}VG_{\phi}^{p-1}\dx<0~\mbox{or}~V\leq 0~\mbox{in}~\Omega.$$
Then~$Q_{V/c_{p}}$ has an optimal Hardy-weight in~$\Omega$.
\ecorollary
The following corollary is adapted from \cite[Corollary~1.2]{Versano}. The proof is similar.
\bcorollary\label{laco}
Let~$V\in\mathbb{M}^{q}_{\loc}(p;\Omega)\cap\widehat{M}^{q}_{\loc}(p;\Omega)$. Suppose that~$Q$ is subcritical in~$\Omega$, that~$V\leq 0$ in~$\Omega$, and that Assumption~\ref{assup} holds. Let~$K$ be an admissible compact subset of~$\Omega$ and let~$G\in C\left(\Omega\setminus \mathring{K}\right)$ be a positive solution $Q'[u]=g$ in~$\Omega\setminus K$ such that
$$G>0\quad\mbox{on}\quad\partial K,\quad\lim_{x\rightarrow\overline{\infty}}G(x)=0,\quad\mbox{and}\quad\int_{\Omega\setminus K}|V|G^{p-1}\dx<\infty,$$
where~$g\in M^{q}_{\loc}(p;\Omega)$ for $p\neq n$, $g\in L^{\bar{\rho}}_{\loc}(\Omega)$ ($\bar{\rho}>1$) for $p=n$, and $g$ is nonnegative a.e. in $\Omega\setminus K$.
Then~$Q_{V/c_{p}}$ has an optimal Hardy-weight in~$\Omega$.
\ecorollary
\bproof
Take~$x_{0}\in\mathring{K}$ and~$0\leq\phi\in\core\setminus\{0\}$ such that~$\supp\phi\subseteq B_{r}(x_{0})\subseteq K$ for some~$r>0$ and that~$\Omega\setminus\supp\phi$ is a domain. By Lemma~\ref{gpcon}, we may find a Green potential~$G_{\phi}$ of~$Q'$ in~$\Omega$ with density~$\phi$. Since~$G_{\phi}\in\mathcal{M}_{\Omega;\,\supp\phi}$ and~$\supp\phi\subseteq\mathring{K}$, there exists a positive constant~$C$ such that~$G_{\phi}\leq CG(x)$ in~$\Omega\setminus K$. Then
$$\lim_{x\rightarrow\overline{\infty}}G_{\phi}(x)=0\quad\mbox{and}\quad\int_{\Omega}|V|G_{\phi}^{p-1}\dx<\infty.$$ By Corollary~\ref{reph},~$Q_{V/c_{p}}$ has an optimal Hardy-weight in~$\Omega$.
\eproof
We conclude the section with two examples realizing the conditions in Corollary~\ref{laco} so that~$Q_{V/c_{p}}$ has an optimal Hardy-weight. See also \cite[Corollary~1.2 and Remark~3.11]{Versano}.
\begin{example}
\emph{Let~$V\in\mathbb{M}^{q}_{\loc}(p;\omega)\cap\widehat{M}^{q}_{\loc}(p;\Omega)$, where~$\omega$ is a Lipschitz domain with~$\omega\Subset\Omega$. Suppose that~$Q$ is nonnegative in~$\Omega$. Then~$Q$ is subcritical in~$\omega$ by Theorem~\ref{perturbation}~(i). By Theorem~\ref{perturbation}~(ii),~$\lambda_{1}(Q;\omega)>0$ (see Definition~\ref{gpeigen}). By Theorem~\ref{peigen}, there exists a positive solution~$G\in W^{1,p}_{0}(\omega)$ of~$Q'[u]=1$ in~$\omega$. By Remark~\ref{nrem},~$G$ is continuous on~$\overline{\omega}$. Then
$$\lim_{x\rightarrow\overline{\infty}}G(x)=0\quad\mbox{and}\quad\int_{\omega}|V|G^{p-1}\dx<\infty,$$
where the~$\overline{\infty}$ is the ideal point in the one-point compactification of~$\omega$. Therefore, if, in addition, in~$\omega$,~$V\leq 0$ and Assumption~\ref{assup} holds, then~$Q_{V/c_{p}}$ has an optimal Hardy-weight in~$\omega$ by Corollary~\ref{laco}.}
\end{example}
\begin{example}
\emph{Let~$\Omega=\R^{n}$, let~$V\in\mathbb{M}^{q}_{\loc}(p;\R^{n})\cap\widehat{M}^{q}_{\loc}(p;\R^{n})$, and let~$1<p<n$. Suppose that the norm family~$H(x,\cdot)$ does not depend on~$x\in\R^{n}$ and that~$\supp V$ is compact. For all~$x\in\R^{n}$, let$$H_{0}(x)=\sup_{\xi\in\R^{n}\setminus\{0\}}\frac{x\cdot\xi}{H(\xi)}$$ be the dual norm of~$H$ and let~$\mathcal{G}(x)=H_{0}(x)^{(p-n)/(p-1)}$ ($x\neq 0$). Then~$\mathcal{G}$ is Finsler $p$-harmonic in~$\R^{n}\setminus\{0\}$ by Examples~\ref{fexas}~(iii). Let~$K$ be an adimissble compact subset of~$\R^{n}$ such that~$\supp V\subseteq K$ and that~$0\in\mathring{K}$. Then$$Q'[\mathcal{G}]=0~\mbox{in}~\R^{n}\setminus K,\quad\mathcal{G}>0~\mbox{on}~\partial K,\quad\lim_{x\rightarrow\infty}\mathcal{G}(x)=0,\quad\mbox{and}\quad\int_{\R^{n}\setminus K}|V|\mathcal{G}^{p-1}\dx=0<\infty.$$ Suppose further that~$Q$ is subcritical, that~$V\leq 0$, and that Assumption~\ref{assup} holds, all in~$\R^{n}$. Then by Corollary~\ref{laco},~$Q_{V/c_{p}}$ has an optimal Hardy-weight in~$\R^{n}$.}

\emph{Moreover, it follows immediately from \cite[Theorem~4.6]{Giri1} that~$Q_{0}$ is subcritical in~$\R^{n}$. Therefore, there exists a nonnegative function~$V_{0}\in L^{\infty}_{c}(\R^{n})\setminus\{0\}$ such that~$Q_{-V_{0}}$ is nonnegative in~$\R^{n}$ (cf.~\cite[Definition~6.3]{HPR}). Furthermore,~$V=-V_{0}/2$ realizes our requirements on the potential term in this example.}
\end{example}
\appendix
\section{Eigenvalues of~$Q'$}\label{eigenvaluea}
Throughout the appendix, we need neither Assumption~\ref{C1alpha} nor~$V\in\widehat{M}^{q}_{\loc}(p;\Omega)$. In other words, we only assume that $(H(x,\cdot))_{x\in\Gw}$ is a family of norms on $\mathbb{R}^{n}$ satisfying Assumptions~\ref{ass9} and that~$V\in M^{q}_\loc(p;\Omega)$.

We first review some definitions of eigenvalues and eigenfunctions of~$Q'$ and two relevant fundamental results, where the second result is used in Sections~\ref{zeropcons} and~\ref{nzpotential}. After some further preliminaries, we prove that the principal eigenvalue of~$Q'$ is isolated and that all the eigenvalues of~$Q'$ form a closed set.
\begin{Def}\label{gpeigen}
	{\em
		The \emph{generalized principal eigenvalue} of $Q'$ in $\Gw$ is defined by		$$\lambda_{1}=\lambda_{1}(Q;\Omega)\triangleq\inf_{\phi\in C^{\infty}_{c}(\Omega) \setminus\{0\}}\frac{Q[\phi;\Omega]}{\Vert \phi\Vert_{L^{p}(\Omega)}^{p}}.$$}	
\end{Def}
\bremark[{\cite[Remark 4.15]{HPR}}]\label{lambda1rem}
\emph{For every domain~$\omega\Subset\Omega$,
$$\lambda_{1}=\inf_{u\in W^{1,p}_{0}(\omega)\setminus\{0\}}\frac{Q[u;\omega]}{\Vert u\Vert_{L^{p}(\omega)}^{p}}.$$}
\eremark
\begin{Def}\label{eigenproblem}
	{\em A real number~$\lambda$ is called an \emph{eigenvalue with an eigenfunction v} of~$Q'$ in~$\Omega$ if~$v\in W^{1,p}_{0}(\Omega)\setminus\{0\}$ is a solution of~$Q'_{V-\lambda}[u]=0$ in~$\Omega$. Such a pair~$(\lambda,v)$ is called an \emph{eigenpair} of~$Q'$ in~$\Omega$. An eigenpair~$(\lambda,v)$ is called \emph{normalized} in $\Omega$  if~$\Vert v\Vert_{L^{p}(\Omega)}=1$.
    }
\end{Def}
\begin{Def}
	{\em  A \emph{principal eigenvalue} of~$Q'$ in~$\Omega$ is an eigenvalue with a nonnegative eigenfunction, which is called  a \emph{principal eigenfunction}.}
\end{Def}
\begin{theorem}[{\cite[Theorem~4.19 and Corollary~5.4]{HPR}}]\label{uniquep}
	In a domain~$\omega\!\Subset\!\Omega$, $\lambda_{1}(Q;\omega)$ is the unique principal eigenvalue of~$Q'$.
\end{theorem}
\begin{theorem}[{\cite[Theorem~4.22 (3) $\Rightarrow$ (5)]{HPR}}]\label{peigen}
Let~$\omega\Subset\Omega$ be a Lipschitz domain. If $\lambda_{1}(Q;\omega)>0$, then for each nonnegative~$g\in L^{p'}(\omega)$, there exists a nonnegative solution~$v\in W^{1,p}_{0}(\omega)$ of~$Q'[u]=g$ in $\omega$ which is either zero (implying~$g=0$) or positive.
\end{theorem}
Next we are going to prove the isolation of~$\lambda_{1}$. In the proof, we will use Corollary~\ref{kato} which is derived from the following result.
\blemma[{\cite[Theorem~2.1]{Ambrosio}}]\label{A21}
If~$v\in W^{1,p}_{\loc}(\Omega)$ is a subsolution of~$Q'[u]=0$ in~$\Omega$, then~$v^{+}$ is a subsolution in~$\Omega$.
\elemma
\bcorollary\label{kato}
If $v\in W^{1,p}_{\loc}(\Omega)$ is a supersolution of~$Q'[u]=0$ in~$\Omega$, then $v^{-}$ is a subsolution in~$\Omega$.
\ecorollary
\begin{proof}
	Note that~$(-v)^{+}=v^{-}$ and then we use Lemma~\ref{A21}.
\end{proof}
Combining the methods of \cite[Theorem~3.9]{PPAPDE}, \cite[Theorem~5.5]{Lnotes}, and \cite[Lemma~6.4]{Regev}, we establish the isolation of the principal eigenvalue~$\lambda_{1}$, bridging a gap in \cite[Theorem~4.19]{HPR}.
We extract part of the proof as the following separate lemma because it will also be used in Theorem~\ref{eigenclose}. See also the proof of \cite[Theorem~3.9]{PPAPDE}.
\blemma\label{eigencon}
Let $\omega\Subset\Omega$ be a domain. Let~$\{(\tilde{\lambda}_{k}, v_{k})\}_{k\in\mathbb{N}}$ be a sequence of normalized eigenpairs  of the operator $Q'$ in~$\omega$ such that~$\lim_{k\rightarrow \infty}\tilde{\lambda}_{k}=\lambda\in \mathbb{R}$. Then~$\{v_{k}\}_{k\in\mathbb{N}}$ is bounded in~$W^{1,p}_{0}(\omega)$ and up to a subsequence, it converges to some function $v\in W^{1,p}_{0}(\omega)$ a.e. in~$\omega$, strongly in~$L^{p}(\omega)$, and weakly in~$W^{1,p}_{0}(\omega)$. 
\elemma
\bproof
By the Morrey-Adams theorem \cite[Theorem~2.4~(i)]{PPAPDE}, we obtain 
$$\int_{\omega}|\nabla v_{k}|_{\mathcal{A}}^{p}\dx\leq |\tilde{\lambda}_{k}|+\int_{\omega}|V||v_{k}|^{p}\dx\leq C_{\delta}+\delta\int_{\omega}|\nabla v_{k}|^{p}\dx$$ for all~$0<\delta<1$ and~$k\in\mathbb{N}$. Furthermore, via some sufficiently small~$\delta>0$, we conclude  that~$\{v_{k}\}_{k\in\mathbb{N}}$ is bounded in~$W^{1,p}_{0}(\omega)$. 
By a simplified version \cite[Lemma~5.1]{Lnotes} of the Rellich-Kondrachov theorem, up to a subsequence, for some~$v\in W^{1,p}_{0}(\omega)$,~$\{v_{k}\}_{k\in\mathbb{N}}$ converges  to~$v$ strongly in~$L^{p}(\omega)$ and weakly in~$W^{1,p}_{0}(\omega)$. Furthermore, up to a subsequence, $\{v_{k}\}_{k\in\mathbb{N}}$ converges to~$v$ a.e. in~$\omega$.
\eproof
Now we can show that~$\lambda_{1}$ is isolated.
\begin{theorem}\label{isolation}
	Let $\omega\Subset\Omega$ be a domain. Then the  principal eigenvalue~$\lambda_{1}(Q;\omega)$ is isolated.
\end{theorem}
\bproof
Suppose that~$\lambda_{1}$ is not isolated. Let~$\{(\tilde{\lambda}_{k}, v_{k})\}_{k\in\mathbb{N}}$ be a sequence of normalized eigenpairs of~$Q'$ such that~$\lim_{k\rightarrow \infty}\tilde{\lambda}_{k}=\lambda_{1}$ and~$\tilde{\lambda}_{k}\neq \lambda_{1}$ for all~$k\in\mathbb{N}$. Then by Lemma~\ref{eigencon},~$\{v_{k}\}_{k\in\mathbb{N}}$ is bounded in~$W^{1,p}_{0}(\omega)$ and up to a subsequence,  for some~$v\in W^{1,p}_{0}(\omega)$, $\{v_{k}\}_{k\in\mathbb{N}}$ converges  to~$v$ a.e. in~$\omega$, strongly in~$L^{p}(\omega)$, and weakly in~$W^{1,p}_{0}(\omega)$. We still denote the corresponding subsequences by~$\{v_{k}\}_{k\in\mathbb{N}}$ and~$\{\tilde{\lambda}_{k}\}_{k\in\mathbb{N}}$. 
Then by \cite[Theorem~4.6]{HPR},
we obtain$$Q[v;\omega]\leq \liminf_{k\rightarrow \infty}Q[v_{k};\omega]=\lambda_{1}\Vert v\Vert_{L^{p}(\omega)}^{p}.$$
By Remark \ref{lambda1rem}, we have
$$Q[v;\omega]\geq\lambda_{1}\Vert v\Vert_{L^{p}(\omega)}^{p}.$$
Then$$Q[v;\omega]=\lambda_{1}\Vert v\Vert_{L^{p}(\omega)}^{p}.$$
By a standard variational argument, $(\lambda_{1},v)$ is a normalized eigenpair. By the proof of \cite[Theorem~4.19]{HPR},~$v$ has a fixed sign. We may assume that~$v$ is positive. 
Otherwise, consider~$\{-v_{k}\}_{k\in \mathbb{N}}$. We follow the methods in \cite[Theorem~3.9]{PPAPDE} and \cite[Lemma~6.4]{Regev}. 
According to Corollary~\ref{kato},~$v^{-}_{k}$ is a subsolution of~$Q'_{V-\tilde{\lambda}_{k}}[v]=0$. 
Then invoking the Morrey-Adams theorem \cite[Theorem~2.4~(i)]{PPAPDE}, we may write, for all~$k\in\mathbb{N}$, all~$\delta>0$, and some positive constant~$C=C(n,p,q)$,
\begin{eqnarray*}
	\int_{\omega}|\nabla v_{k}^{-}|_{\mathcal{A}}^{p}\dx&\leq&\int_{\omega}|V-\tilde{\lambda}_{k}||v_{k}^{-}|^{p}\dx\\
	&\leq& \delta\Vert \nabla v_{k}^{-}\Vert^{p}_{L^{p}(\omega;\,\mathbb{R}^{n})}+C\delta^{-n/(pq-n)}\Vert V-\tilde{\lambda}_{k}\Vert_{M^{q}(p;\,\omega)}^{pq/(pq-n)}\Vert v_{k}^{-}\Vert^{p}_{L^{p}(\omega)}.
\end{eqnarray*}
For all~$0<\delta<\alpha_{\omega}$ and~$k\in\mathbb{N}$, by virtue of the local uniform ellipticity in Theorem~\ref{thm_1}, we get
$$(\alpha_{\omega}-\delta)\Vert \nabla v_{k}^{-}\Vert^{p}_{L^{p}(\omega;\,\mathbb{R}^{n})}\leq C\delta^{-n/(pq-n)}\Vert V-\tilde{\lambda}_{k}\Vert_{M^{q}(p;\,\omega)}^{pq/(pq-n)}\Vert v_{k}^{-}\Vert^{p}_{L^{p}(\omega)}.$$
By H\"older's inequality, for all~$0<\delta<\alpha_{\omega}$ and~$k\in\mathbb{N}$,
$$(\alpha_{\omega}-\delta)\Vert \nabla v_{k}^{-}\Vert^{p}_{L^{p}(\omega;\,\mathbb{R}^{n})}\leq C\delta^{-n/(pq-n)}\Vert V-\tilde{\lambda}_{k}\Vert_{M^{q}(p;\,\omega)}^{pq/(pq-n)}\vol(v_{k}^{-1}((-\infty,0)))^{1-p/\bar{p}}\Vert v_{k}^{-}\Vert^{p}_{L^{\bar{p}}(\omega)},$$
where~$\bar{p}\triangleq
\begin{cases}
	2p&\mbox{if}~p\geq n,\\
	\frac{np}{n-p}&\mbox{if}~p<n.
\end{cases}$ 
By \cite[Corollary~3.7]{Kinnunen}, for all~$0<\delta<\alpha_{\omega}$ and~$k\in\mathbb{N}$, 
$$(\alpha_{\omega}-\delta)\Vert \nabla v_{k}^{-}\Vert^{p}_{L^{p}(\omega;\,\mathbb{R}^{n})}\leq C\delta^{-n/(pq-n)}\Vert V-\tilde{\lambda}_{k}\Vert_{M^{q}(p;\,\omega)}^{pq/(pq-n)}\vol(v_{k}^{-1}((-\infty,0)))^{1-p/\bar{p}}\Vert \nabla v_{k}^{-}\Vert^{p}_{L^{p}(\omega)},$$ where the constant~$C$ is independent of~$k\in\mathbb{N}$. By Theorem~\ref{uniquep}, for all~$k\in\mathbb{N}$,~$v^{-}_{k}\neq 0$ and hence~$\Vert\nabla v_{k}^{-}\Vert^{p}_{L^{p}(\omega)}\neq 0.$ Then for all~$0<\delta<\alpha_{\omega}$ and~$k\in\mathbb{N}$,
$$(\alpha_{\omega}-\delta)\leq C\delta^{-n/(pq-n)}\Vert V-\tilde{\lambda}_{k}\Vert_{M^{q}(p;\,\omega)}^{pq/(pq-n)}\vol(v_{k}^{-1}((-\infty,0)))^{1-p/\bar{p}}.$$
With the help of the positive principal eigenfunction and Poincar\'e's inequality, we must have~$\Vert V-\lambda_{1}\Vert_{M^{q}(p;\,\omega)}>0$, which implies that  for all sufficiently large~$k\in\mathbb{N}$,~$\Vert V-\tilde{\lambda}_{k}\Vert_{M^{q}(p;\,\omega)}>0.$ Now a positive lower bound of~$\vol(v_{k}^{-1}((-\infty,0)))$ for all sufficiently large~$k\in\mathbb{N}$ has been clear, which prohibits~$\lim_{k\rightarrow\infty}\vol(v_{k}^{-1}((-\infty,0)))=0.$

The rest is similar to the counterpart in \cite[Theorem~3.9]{PPAPDE}. Take an arbitrary~$\varepsilon>0$ and a compact subset~$K$ of~$\omega$ such that~$\vol(\omega\setminus K)<\varepsilon$. Since~$v$ is continuous and positive in~$\omega$,~$v$ has a positive lower bound on~$K$. Because~$\{v_{k}\}_{k\in\mathbb{N}}$ converges to~$v$ a.e. in~$\omega$, by \cite[Theorem~11.32]{Hewitt}, there exists a measurable set~$K'\subseteq K$ such that~$\vol(K')<\varepsilon$ and~$\{v_{k}\}_{k\in\mathbb{N}}$ converges uniformly to~$v$ on~$K\setminus K'$. Now that~$v$ has a positive lower bound on~$K$, we may conclude that for all sufficiently large~$k\in\mathbb{N}$, the functions~$v_{k}$ are nonnegative on~$K\setminus K'$. Hence for all sufficiently large~$k\in\mathbb{N}$,~$v_{k}^{-1}((-\infty,0))\subseteq K'\cup (\omega\setminus K).$ Furthermore, for all sufficiently large~$k\in\mathbb{N}$,~$\vol(v_{k}^{-1}((-\infty,0)))<2\varepsilon$. Then~$\lim_{k\rightarrow\infty}\vol(v_{k}^{-1}((-\infty,0)))=0.$  Contradiction! 
\eproof
The following theorem will be used in Theorem~\ref{eigenclose}. Recall that $1<p<\infty$.
\begin{theorem}[{\cite[Theorem~13.44]{Hewitt}}]\label{wconvergence}
Let $(X,\mathscr{A},\mu)$ be a measure space, where $\mathscr{A}$ is a $\sigma$-algebra of subsets of the set $X$ and $\mu$ is a measure on $\mathscr{A}$. Suppose that $\{f_{k}\}_{k\in\mathbb{N}}$ is a bounded sequence in $L^{p}(X,\mathscr{A},\mu)$ and that $f_{k}\xrightarrow{k \to \infty}f$ $\mu$-a.e.. Then $f_{k}\xrightarrow{k \to \infty} f$ weakly in $L^{p}(X,\mathscr{A},\mu)$.
\end{theorem}
Inspired by the methods of \cite[Theorem~5.1]{Lnotes} and \cite[Theorem~3.9]{PPAPDE}, we prove:
\btheorem\label{eigenclose}
All the real eigenvalues of~$Q'$ in a domain~$\omega\Subset\Omega$ form a closed subset of~$\R$.
\etheorem
\bproof
Let~$\{(\tilde{\lambda}_{k}, v_{k})\}_{k\in\mathbb{N}}$ be a sequence of normalized eigenpairs of~$Q'$ such that$$\lim_{k\rightarrow \infty}\tilde{\lambda}_{k}=\lambda\in\R.$$ By Lemma~\ref{eigencon}, $\{v_{k}\}_{k\in\mathbb{N}}$ is bounded in~$W^{1,p}_{0}(\omega)$ and up to a subsequence, $\{v_{k}\}_{k\in\mathbb{N}}$ converges to a function $v\in W^{1,p}_{0}(\omega)$ a.e. in~$\omega$, strongly in~$L^{p}(\omega)$, and weakly in~$W^{1,p}_{0}(\omega)$. We continue writing the corresponding subsequences as~$\{v_{k}\}_{k\in\mathbb{N}}$ and~$\{\tilde{\lambda}_{k}\}_{k\in\mathbb{N}}$.
Note that by \cite[Remark~6.18]{Hou},
\begin{eqnarray*}
	&&\int_{\omega}\big(\mathcal{A}(x,\nabla v_{k})-\mathcal{A}(x,\nabla v)\big)\cdot \nabla(v_{k}-v)\dx+\int_{\omega}V|v_{k}|^{p-2}v_{k}(v_{k}-v)\dx\\
	&=&\tilde{\lambda}_{k}\int_{\omega}|v_{k}|^{p-2}v_{k}(v_{k}-v)\dx-\int_{\omega}\mathcal{A}(x,\nabla v)\cdot \nabla(v_{k}-v)\dx.
\end{eqnarray*}
Because $\{v_{k}\}_{k\in\mathbb{N}}$ converges to~$v$ weakly in~$W^{1,p}_{0}(\omega)$ and strongly in~$L^{p}(\omega)$, the right hand side converges to zero. By the H\"older inequality, we get $$\int_{\omega}|V||v_{k}|^{p-1}|v_{k}-v|\dx
\leq\left(\int_{\omega}|V||v_{k}|^{p}\dx\right)^{1-1/p}\left(\int_{\omega}|V||v_{k}-v|^{p}\dx\right)^{1/p}.$$ Since~$\{v_{k}\}_{k\in\mathbb{N}}$ is bounded in~$W^{1,p}_{0}(\omega)$, by the Morrey-Adams theorem \cite[Theorem~2.4~(i)]{PPAPDE}, $\left(\int_{\omega}|V||v_{k}|^{p}\dx\right)^{1-1/p}$ is bounded with respect to~$k\in\mathbb{N}$. For each~$\delta>0$ and all sufficiently large~$k\in\mathbb{N}$, by the Morrey-Adams theorem \cite[Theorem~2.4~(i)]{PPAPDE}, we obtain
\begin{eqnarray*}
	\int_{\omega}|V||v_{k}-v|^{p}\dx&\leq& \delta\int_{\omega}|\nabla(v_{k}-v)|^{p}\dx+C\delta^{-n/(pq-n)}\int_{\omega}|v_{k}-v|^{p}\dx\\	&\leq&\delta\int_{\omega}|\nabla(v_{k}-v)|^{p}\dx+C\delta.
\end{eqnarray*}
Because~$\{v_{k}\}_{k\in\mathbb{N}}$ is bounded in~$W^{1,p}_{0}(\omega)$, we deduce that~$\lim_{k\rightarrow\infty}\int_{\omega}|V||v_{k}-v|^{p}\dx=0.$ Hence,
\begin{align*}
\lim_{k\rightarrow\infty}\int_{\omega}\left(\mathcal{A}(x,\nabla v_{k})-\mathcal{A}(x,\nabla v)\right)\cdot\nabla(v_{k}-v)\dx=\lim_{k\rightarrow\infty}\int_{\omega}V|v_{k}|^{p-2}v_{k}(v_{k}-v)\dx=0.
\end{align*}
Using \cite[Lemma~3.73]{HKM}, we deduce that~$\{\mathcal{A}(x,\nabla v_{k})\}_{k\in\mathbb{N}}$ converges weakly to~$\mathcal{A}(x,\nabla v)$ in~$L^{p'}(\omega;\R^{n})$. Then for all~$\phi\in C^{\infty}_{c}(\omega)$,
\begin{align}\label{convergence1}
\lim_{k\rightarrow\infty}\int_{\omega}\mathcal{A}(x,\nabla v_{k})\cdot\nabla\phi\dx=\int_{\omega}\mathcal{A}(x,\nabla v)\cdot\nabla\phi\dx.
\end{align}
It can be easily seen that~$\{|v_{k}|^{p-2}v_{k}\}_{k\in\mathbb{N}}$ is bounded in~$L^{p'}(\omega)$ because~$\{v_{k}\}_{k\in\mathbb{N}}$ is bounded in~$L^{p}(\omega)$.
Since~$\{v_{k}\}_{k\in\mathbb{N}}$ converges to~$v$ a.e. in~$\omega$,~$\{|v_{k}|^{p-2}v_{k}\}_{k\in\mathbb{N}}$ converges to~$|v|^{p-2}v$ a.e. in~$\omega$. By Theorem~\ref{wconvergence},~$\{|v_{k}|^{p-2}v_{k}\}_{k\in\mathbb{N}}$ converges weakly to~$|v|^{p-2}v$ in~$L^{p'}(\omega)$. Then due to~$\lim_{k\rightarrow \infty}\tilde{\lambda}_{k}=\lambda\in\R$, for all~$\phi\in C^{\infty}_{c}(\omega)$,
\begin{align}\label{convergence2}
\lim_{k\rightarrow\infty}\tilde{\lambda}_{k}\int_{\omega}|v_{k}|^{p-2}v_{k}\phi\dx=\lambda\int_{\omega}|v|^{p-2}v\phi\dx.
\end{align}
By the Morrey-Adams theorem \cite[Theorem~2.4~(i)]{PPAPDE},
it is a simple matter to check that~$\{\int_{\omega}V^{\pm}|v_{k}|^{p}\dx\}_{k\in\mathbb{N}}$ are bounded because~$\{v_{k}\}_{k\in\mathbb{N}}$ is bounded in~$W^{1,p}_{0}(\omega)$. Similarly, since~$\{|v_{k}|^{p-2}v_{k}\}_{k\in\mathbb{N}}$ converges to~$|v|^{p-2}v$ a.e. in~$\omega$, by Theorem~\ref{wconvergence}, for all~$\phi\in C^{\infty}_{c}(\omega)$,
\begin{align*}
\lim_{k\rightarrow\infty}\int_{\omega}V^{\pm}|v_{k}|^{p-2}v_{k}\phi\dx=\int_{\omega}V^{\pm}|v|^{p-2}v\phi\dx.
\end{align*}
Then for all~$\phi\in C^{\infty}_{c}(\omega)$,
\begin{align}\label{convergence3}
\lim_{k\rightarrow\infty}\int_{\omega}V|v_{k}|^{p-2}v_{k}\phi\dx=\int_{\omega}V|v|^{p-2}v\phi\dx.
\end{align}
Recall that for all~$k\in\mathbb{N}$ and~$\phi\in C^{\infty}_{c}(\omega)$,
\begin{align}\label{convergence4}
\int_{\omega}\mathcal{A}(x,\nabla v_{k})\cdot\nabla\phi\dx+\int_{\omega}V|v_{k}|^{p-2}v_{k}\phi\dx=\tilde{\lambda}_{k}\int_{\omega}|v_{k}|^{p-2}v_{k}\phi\dx.
\end{align}
Combining \eqref{convergence1}, \eqref{convergence2}, \eqref{convergence3}, and \eqref{convergence4}, we deduce that for all~$\phi\in C^{\infty}_{c}(\omega)$,
\begin{align*}
\int_{\omega}\mathcal{A}(x,\nabla v)\cdot\nabla\phi\dx+\int_{\omega}V|v|^{p-2}v\phi\dx=\lambda\int_{\omega}|v|^{p-2}v\phi\dx.
\end{align*}
Since~$\Vert v_{k}\Vert_{L^{p}(\omega)}=1$ for all~$k\in\mathbb{N}$ and $\{v_{k}\}_{k\in\mathbb{N}}$ converges to~$v$ in~$L^{p}(\omega)$, we conclude that~$\Vert v\Vert_{L^{p}(\omega)}=1$. Then~$v\not\equiv 0$ in~$\omega$. Therefore,~$\lambda$ is an eigenvalue of~$Q'$ in~$\omega$ and the desired conclusion follows.
\eproof
\section*{Acknowledgments}
This work is part of the author’s ongoing Ph.D. thesis of the Technion - Israel Institute of Technology written under the supervision of Professors Yehuda Pinchover and Matthias Keller. The author extends deep gratitude to Professor Pinchover for his invaluable guidance throughout the preparation of this paper and to Prof. Dr. Keller for his advice on the appendix. The author also acknowledges the financial support provided by the Technion and the Israel Science Foundation (Grant No. 637/19) founded by the Israel Academy of Sciences and Humanities.
{}
\end{document}